%% file: main.tex
\documentclass[preprint,sort&compress,number]{elsarticle}

\usepackage{amsmath,amsthm,amssymb,latexsym,amscd}
\usepackage{subfig}
\usepackage[list-final-separator={, and }]{siunitx}
\usepackage[shadow,textsize=small,textwidth=2.5cm]{todonotes}
\usepackage{stmaryrd}

\usepackage{geometry}
\usepackage{mathtools}
\usepackage{algorithm}
\usepackage{algpseudocode}
\algtext*{EndIf}
\algtext*{EndWhile}
\algtext*{EndFor}

\usepackage{color}
\usepackage{nomencl}
\usepackage{bm}
\usepackage{lineno}
\usepackage{tablefootnote}
\usepackage{longtable}
\usepackage{threeparttable}
\usepackage{caption}
\makeatletter 
\def\@makecaption#1#2{%
  \vskip\abovecaptionskip
  \sbox\@tempboxa{#1 #2}%
    {\bfseries #1} #2\par
  \vskip\belowcaptionskip}
\makeatother
\usepackage{booktabs}
\usepackage[justification=centering]{caption}
\usepackage{setspace}

\usepackage{amsmath,amsthm,amssymb,latexsym,amscd}
\usepackage{lineno}

\usepackage{changes}
\usepackage{verbatim}
\usepackage{todonotes}

\input{mycommand.tex}
\usepackage{hyperref}
\hypersetup{
	colorlinks=true,
	linkcolor=blue,
	filecolor=magenta,      
	urlcolor=cyan,
}

    \title{A phase-field fracture model in thermo-poro-elastic media with micromechanical strain energy degradation}
\author[TON1,TON2]{Yuhao Liu}
\author[MUL,UFZ]{Keita Yoshioka}
\author[MUL,UFZ]{Tao You}
\author[TON1,TON2]{Hanzhang Li}
\author[TON1,TON2]{Fengshou Zhang}

\address[TON1]{Key Laboratory of Geotechnical $\&$ Underground Engineering of Ministry of Education, Tongji University, Shanghai 200092, China}
\address[TON2]{Department of Geotechnical Engineering, College of Civil Engineering, Tongji University, Shanghai 200092, China}
\address[MUL]{Department Geoenergy, Montanuniversit\"at Leoben, Leoben 8700, Austria}
\address[UFZ]{Department of Environmental Informatics, Helmholtz Centre for Environmental Research - UFZ, Leipzig 04318, Germany}

\begin{document}
\captionsetup[figure]{labelfont={bf},name={Fig.},labelsep=period}
\begin{abstract}
This work extends the hydro-mechanical phase-field fracture model to non-isothermal conditions with micromechanics based poroelasticity, which degrades Biot’s coefficient not only with the phase-field variable (damage) but also with the energy decomposition scheme. 
Furthermore, we propose a new approach to update porosity solely determined by the strain change rather than damage evolution as in the existing models. 
As such, these poroelastic behaviors of Biot’s coefficient and the porosity dictate Biot’s modulus and the thermal expansion coefficient. 
For numerical implementation, we employ an isotropic diffusion method to stabilize the advection-dominated heat flux and adapt the fixed stress split method to account for the thermal stress.
We verify our model against a series of analytical solutions such as Terzaghi's consolidation, thermal consolidation, and the plane strain hydraulic fracture propagation, known as the KGD fracture. 
Finally, numerical experiments demonstrate the effectiveness of the stabilization method and intricate thermo-hydro-mechanical interactions during hydraulic fracturing with and without a pre-existing weak interface. 
\vspace{2em}

\end{abstract}

\maketitle

\section{Introduction}
\input{intro}

\section{Mathematical model}
\label{sec:mat_model}
\input{mathe}

\section{Numerical implementation}
\label{sec:num_imp}
\input{numeric}

\section{Model verification}
\label{sec:verification}
\input{verif}

\section{Numerical experiments and discussion}
\label{sec:num_exp}
\input{numexperi}

\section{Conclusions}
\label{sec:concl}
\input{conc}
\section*{Declaration of competing interest}
The authors declare that they have no known competing financial interests or personal relationships that could
have appeared to influence the work reported in this paper.

\section*{Data availability}
No data was used for the research described in the article.

\section*{Acknowledgement}
This work is supported by the National Key Research and Development Project (No.2023YFE0110900) and National Natural Science Foundation of China (No.42320104003, 42077247).

\section*{Appendix}
\input{appendix}
\bibliographystyle{apalike}
\bibliography{references} 
\end{document}

%% file: mycommand.tex
\usepackage{amsmath,amsthm,amssymb,latexsym,amscd}
\usepackage{algorithm}
\usepackage{algpseudocode}
\usepackage{lineno}

\newcommand{\ATone}{\texttt{AT$_1$}}
\newcommand{\ATtwo}{\texttt{AT$_2$}}

\newcommand{\eps}{\varepsilon}

\newcommand{\mathd}{\mathrm{d}}
\newcommand{\mrm}{\mathrm}
\newcommand{\mbf}[1]{{\mathbf{#1}}}

\newcommand{\mbfs}{\boldsymbol}
\newcommand{\Tr}[1]{\mrm{Tr} \left( #1\right) }

\DeclareMathOperator*{\argmin}{argmin}

%% file: intro.tex

Nucleation and propagation of fractures are one of the most significant interests or concerns for geo-energy applications such as geothermal/hydrocarbon production, energy storage, or CO$_2$/nuclear waste disposal~\cite{detournay1995coupled, li2006thermo, pandey2018geothermal, lee2020numerical}.
As we exploit deeper subsurface formations, higher surrounding pressures and temperatures impose higher differential pressures and temperatures in the multi-physical processes of fracturing~\cite{hamza2023introduction}.
In deep strata, one may encounter a lower critical pressure~\cite{zeng2020extended} that can lead to larger crack openings~\cite{luo2022numerical, feng2016numerical}, and at high temperatures, thermal branching fractures may initiate normal to the main fracture~\cite{jiao2022investigation, yan2022fdem}. 
Also, thermal erosion may cause cracks to penetrate the caprock in CO$_2$ sequestration operations~\cite{fu2021thermo}. 
Furthermore, the thermal stress can result in tensile fracture following the cooled area for long-term water injection~\cite{parisio2019risks}. 

Because of high pressure and temperature conditions and the large scale involved in geoenergy applications, numerical modeling has been instrumental in forecasting fracture evolution and analyzing seal integrity over time~\cite{bakhshi2021numerical}.
In the numerical simulations of thermo-hydro-mechanical hydraulic fracture, two main challenges still exist. 
One is to properly capture the interactions between fluid flow, heat transfer and fracture mechanics -- particularly the fracturing criteria~\cite{chaudhuri2009buoyant, wang2003coupled}.
The other challenge is the numerical representation of fractures, sharp discontinuities that evolve~\cite{chukwudozie2019variational}. 
Fractures often propagate in a complex topology in rocks, but numerous models impose limitations on the propagation paths. 
Cohesive zone models~\cite{carrier2012numerical} insert lower-dimensional interface elements along the crack surface, constraining the fracture propagation to the prescribed surfaces.
Extended finite element methods~\cite{belytschko1999elastic, lecampion2009extended} address this drawback but struggle to simulate fracture branching and merging. 
Discrete element~\cite{potyondy2004bonded} and lattice based methods~\cite{damjanac2016application} still pose the weaknesses on element dependency, parameter calibration, and calculation scale.

Phase-field models of fracture~\cite{bourdin2000numerical, bourdin2008variational} have emerged as an effective computational approach in the last couple of decades, overcoming many of these shortcomings.
For fluid-driven (hydraulic) fracture, phase-field modeling has been first applied by Bourdin~et~al. (2012)~\cite{bourdin2012variational} and Chukwudozie~et~al. (2013)~\cite{chukwudozie2013variational} for elastic media and then further extended to poroelastic media by Wheeler et al. (2014)~\cite{wheeler2014augmented} and Mikelic~et~al. (2015)~\cite{mikelic2015quasi}. 
Subsequently, many studies have been carried out on the hydro-mechanical phase-field method over various aspects, including non-porous fluid flows~\cite{miehe2015minimization, wilson2016phase}, crack opening calculation~\cite{lee2017iterative, yoshioka2020crack} and approximate treatment for hydraulic parameters~\cite{miehe2016phase, xia2017phase}. 

For thermo-hydro-mechanical (THM) coupled modeling, which is the main topic in this study, Noii and Wick (2019)~\cite{noii2019phase} first extended the hydro-mechanical phase-field fracture model to include thermal effects by introducing the work of temperature from thermo-poroelasticity in the total energy functional.
However, their model neglected spatial variations of the pressure field and, consequently, the convective heat transfer.
For a more complete coupling of thermo-hydro-mechanical processes, Li~et~al.~\cite{li2021phase} implemented the multi-physics interactions during the fracture propagation. 
Suh and Sun (2021)~\cite{suh2021asynchronous} proposed a thermo-hydro-mechanical phase-field model and applied an asynchronous operator-split framework to capture the heat transfer between the fracture and matrix. 
Yi et al. (2024)~\cite{yi2024coupled} conducted numerical experiments of thermo-hydro-mechanical hydraulic fractures under more complicated settings like natural fracture networks, temperature variations, and multi-mineral rocks. 
In addition to the thermo-hydro-mechanical coupling, Feng et al. (2023)~\cite{feng2023phase} introduced effects of CO$_2$ fracturing, and Dai et al. (2024)~\cite{dai2024thermal} incorporated reactive chemistry to simulate acid fracturing.

Despite all the recent developments in thermo-hydro-mechanical coupled phase-field fracture models, the existing works have not adequately addressed the following three aspects.
The first is the degradation of poro-thermo-elastic strain energy.
While the degradation of elastic strain energy has been more or less established, somewhat ad-hoc and various ways to degrade the poro-thermo-elastic strain energy have been employed in~\cite{noii2019phase, suh2021asynchronous, wang2023framework, yi2024coupled}.
The second is the transition of the hydraulic and thermal material parameters from intact to fractured material.
Although the most popular approach is to use a linear indicator function~\cite{li2021phase, feng2023phase}, a linear interpolation does not conserve the mass or energy.
Lastly, the advection-dominated heat transfer in fractures can cause numerical oscillation~\cite{wang2003coupled}, and some stabilization methods may be required~\cite{wang2023framework}. 


In this study, we aim to address these three aspects.
First, we extend the existing hydro-mechanical phase-field fracture model to thermo-poro-elastic media by introducing the thermal strain.
The degradation of the thermo-poro-elastic strain energy is achieved through micromechanically derived Biot's coefficient, and its degradation depends on the phase-field (damage) and the energy decomposition scheme. 
Then, we propose a new way to update the porosity that only depends on the strain change rather than the damage.
From the degraded poro-elastic coefficients and porosity, we approximate hydraulic and thermal properties that transition from a fully-damaged to an intact state of the material.
Finally, we employ an isotropic diffusion method to stabilize the advective dominated heat flow in fracture.
For numerical implementation, we apply a hybrid staggered scheme and adapt the fixed stress split stabilization method by accounting for thermal stress. 
Our proposed model is then tested by a series of analytical solutions verifying different implementation modules. 

The structure of this paper is as follows. 
In Section~\ref{sec:mat_model}, we propose a thermo-hydro-mechanical phase-field fracture model in thermo-poro-elastic media.
Section~\ref{sec:num_imp} discusses the numerical implementation of the new proposed model with the isotropic diffusion method and the adapted fixed stress split strategy for the thermo-hydro-mechanical phase-field coupling problem. 
In section~\ref{sec:verification}, we verify the proposed model against a series of analytical solutions. 
The accuracy of the newly proposed porosity update is highlighted against the existing hydro-mechanical phase-field model. 
Section~\ref{sec:num_exp} studies the impacts of thermo-hydro-mechanical on hydraulic fracture, followed by conclusions and possible future studies.

We use the following notations throughout the paper. The second-order identity tensor is denoted by $\mbf{I}$. 
The symbol $||\mbf{A}|| = \sqrt{\mbf{A}:\mbf{A}}$ is used to calculate the norm of any second-order tensor $\mbf{A}$. The fourth order projectors $\mathbb{J}$ and $\mathbb{K}$ are expressed in the component form as $J_{ijkl}=(\delta_{ik}\delta_{jl}/3$ and $K_{ijkl}=(\delta_{ik}\delta_{jl}+\delta_{il}\delta_{jk})/2-\delta_{ij}\delta_{kl}/3$. The trace operator $\mathrm{Tr}(\cdot)$ acting on the second-order tensors $\mbf{A}$ is defined as $\mathrm{Tr}(\mbf{A}) = \mbfs{\delta}:\mbf{A}$. $\nabla(\cdot)$ is the gradient of $(\cdot)$.



%% file: mathe.tex
\begin{figure}[H]
		\centering
		\includegraphics[scale=0.4]{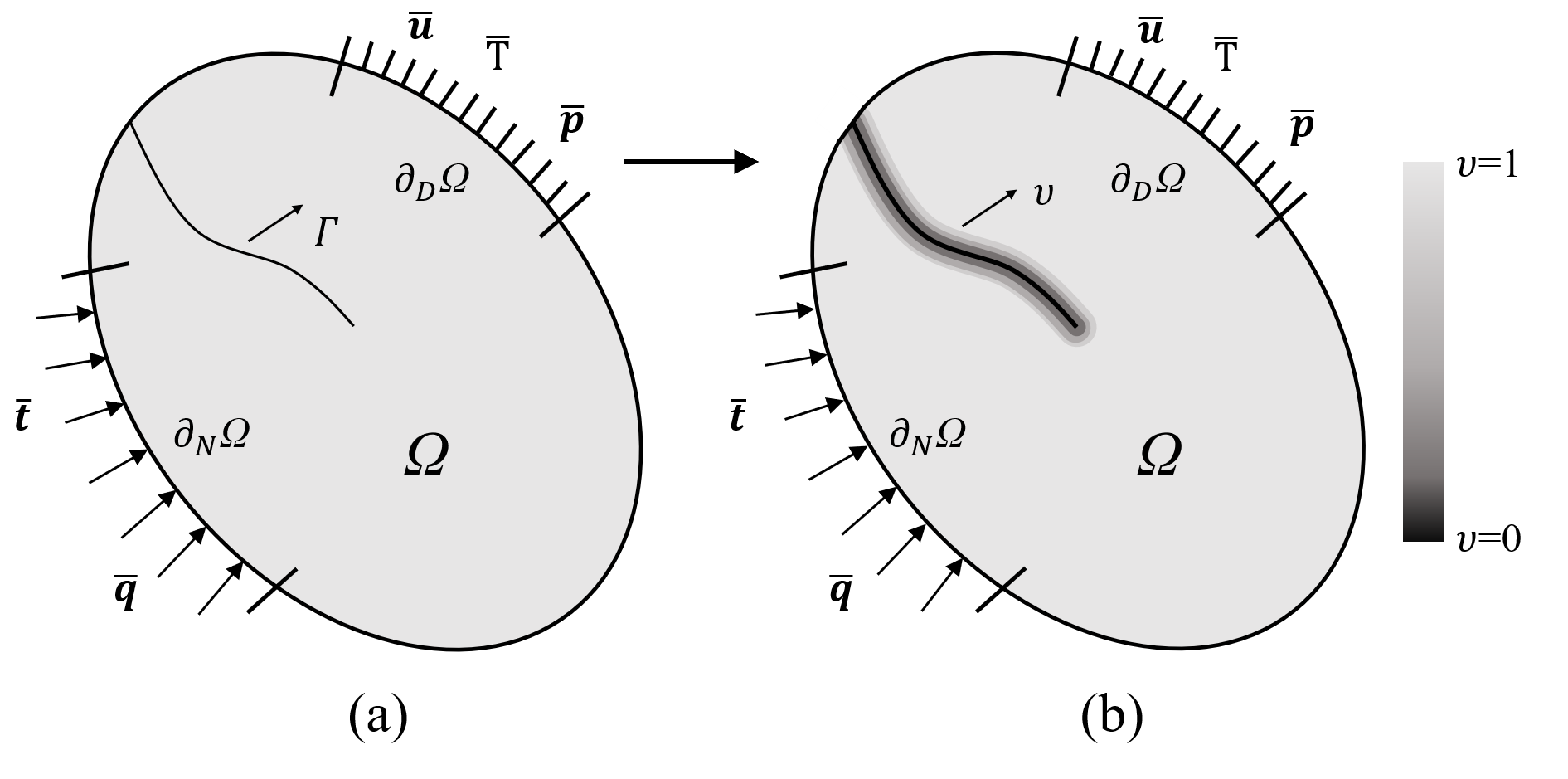}
		\caption{(a) The schematic of porous medium with an existing fracture and (b) the diffused representation of a fracture $\mit\mit{\Gamma}$ in the poroelastic medium.}
		\label{fig:pf_model}
\end{figure}

Consider a thermo-poro-elastic medium that occupies the domain $\mit{\Omega}\in\mathbb{R}^d$, $d=2,3$. In $\mit{\Omega}$, a lower dimensional set of fractures is denoted by $\mit{\Gamma}\in\mathbb{R}^{d-1}$ (Fig.~\ref{fig:pf_model}a). 
The body is subjected to a possible flux $\bar{\mbf{q}}$ and a surface force $\bar{\mbf{t}}$ on the boundary $\mit\partial_\mrm{N}\mit{\Omega} := \mit{\mathcal{C}}_\mrm{N}\cup\partial\mit{\Gamma}$, where $\mathcal{C}_\mrm{N}$ denotes the outer domain boundary and $\partial\mit{\Gamma}$ the fracture boundary. 
The prescribed displacement $\bar{\mbf{u}}$, pressure $\bar{\mbf{p}}$ and temperature $\bar{\mbf{T}}$ can be applied on the boundary $\mit\partial\mit{\Omega}_\mrm{D}$. 
We consider a single-phase Newtonian fluid (liquid) in both pore and fracture spaces, and for injected fluid. 
Furthermore, we assume a local thermal equilibrium in the domain in which the fluid and porous solid temperatures are equilibrated. 

We first derive the mechanical equilibrium in the setting with a discrete crack set $\Gamma$ (Fig.~\ref{fig:pf_model}a).
Then, we regularize the discrete crack set with the now-classical phase-field approach~\cite{bourdin2000numerical, bourdin2008variational} and derive equivalent hydraulic and thermal properties in the setting with diffused cracks (Fig.~\ref{fig:pf_model}b) so that the conservation laws apply over the entire domain $\mit\Omega$.
Subsequently, mass and heat balances are derived with the equivalent properties in the diffused crack setting.

\subsection{Mechanical equilibrium with discrete cracks}
The mechanical equilibrium and fracture evolution in the thermo-poroelastic medium can be derived based on the variational approach proposed by Francfort and Marigo (1998)~\cite{francfort1998revisiting}. 
They defined the total energy of a homogeneous isotropic and linearly elastic material with a set of fractures as a sum of the strain energy, the surface energy of fracture, the body force $\mbf{b}$ and the traction $\mbf{\bar{t}}$ on the boundary $\mit{\mathcal{C}}_\mrm{N}$:
\begin{equation}
    \label{eq:energy_functional}
    \mathcal{E}_\Gamma(\mbf{u}, \mit{\Gamma}) :
    = \int_{\mit{\Omega}\setminus\mit\Gamma}^{}\psi(\mbf{u}) \, \mathd V + \int_{\mit{\Gamma}}^{}G_c \, \mathd S - \int_{\mit{\Omega}\setminus\mit\Gamma}^{}\mbf{b}\cdot\mbf{u} \, \mathd V - \int_{\mit{\mathcal{C}}_\mrm{N}}^{}\mbf{\bar{t}}\cdot\mbf{u} \, \mathd V
\end{equation}
where $G_c$ is the critical surface energy release rate, $\mbf{u}$ is the displacement field, and $\psi(\mbf{u})$ is the strain energy density described by the elastic stiffness tensor $\mathbb{C}_\mrm{m}$ and linearized strain $\mbfs{\eps}(\mbf{u}) = \frac{1}{2} (\nabla \mbf{u} + \nabla \mbf{u}^{\text{T}} )$:
\begin{equation}
    \label{eq:strain_energy}
    \psi(\mbf{u}) = \frac{1}{2} \mathbb{C}_{m} : \mbfs{\eps}(\mbf{u}) : \mbfs{\eps}(\mbf{u})
    .
\end{equation} 

Following~\cite{miehe2015phase, suh2021asynchronous}, we propose the form of total energy in thermo-poroelastic medium
\begin{equation}
    \label{eq:thermo_energy_functional}
    \mathcal{F}_\Gamma(\mbf{u}, p, T, \mit{\Gamma}) := W(\mbf{u}, p, T) + \int_{\mit{\Gamma}}^{}G_c \, \mathd S - \int_{\mit{\Omega}\setminus\mit\Gamma}^{}\mbf{b}\cdot\mbf{u} \, \mathd V - \int_{\mit{\mathcal{C}}_\mrm{N}}^{}\mbf{\bar{t}}\cdot\mbf{u} \, \mathd V
    ,
\end{equation}
where thermo-poro-elastic strain energy $W$ is assumed to decompose into elastic, hydraulic, and thermal parts:
\begin{equation}
        \label{eq:thermo_poroelastic_strain_energy_total}
        W(\mbf{u}, p, T) := W_{\mathrm{elastic}}(\mbf{u}, T) + W_{\mathrm{fluid}}(\mbf{u}, \zeta) + W_{\mathrm{thermo}}(T).
\end{equation}

The thermo-elastic strain energy is generated by the intergranular stress acting on the solid skeleton and is given as:
\begin{equation}
        \label{eq:thermo_elastic_strain_energy}
        W_{\mathrm{elastic}}(\mbf{u}, T) := \int_{\mit{\Omega}\setminus\mit\Gamma}^{}\psi_e(\mbf{u}, T) \, \mathd V.
\end{equation}
Thermo-poro-elastic analyses in a small strain theory show that the total strain is decomposed into the elastic strain $\eps_e$ and the thermal strain $\eps_t$ induced by the temperature difference between the current temperature and the initial temperature $T_0$, i.e., $\Delta T = T - T_0$~\cite{selvadurai2016thermo}:
\begin{equation}
    \label{eq:total_strain}
    \mbfs{\eps} = \mbfs{\eps_e} + \mbfs{\eps_t} = \mbfs{\eps_e} + \beta\Delta T\mbf{I}
    ,
\end{equation}
where $\beta$ is the thermal expansion coefficient.
The elastic strain energy, $ \psi(\mbf{u})$, is solely contributed by the elastic strain~\cite{badnava2018h} and thus we have
\begin{equation}
    \label{eq:thermo_elastic_strain_energy_density}
    \psi_e(\mbf{u}, T) = \frac{1}{2} \mathbb{C}_{m} : \mbfs{\eps}_e(\mbf{u}) : \mbfs{\eps}_e(\mbf{u}) = \frac{1}{2} \mathbb{C}_{m} : \Big(\mbfs{\eps}(\mbf{u})-\beta\Delta T\mbf{I}\Big) : \Big(\mbfs{\eps}(\mbf{u})-\beta\Delta T\mbf{I}\Big).
\end{equation}

Following ~\cite{coussy2004poromechanics,li2021phase, YOU2023116305,yi2024coupled}, the energy of pore fluid can be defined as
\begin{equation}
    \label{eq:pore_fluid_energy}
    W_{\mathrm{fluid}}(\mbf{u}, \zeta) := \int_{\mit{\Omega}\setminus\mit\Gamma}^{}\psi_f(\mbf{u}. \zeta) \, \mathd V = \int_{\mit{\Omega}\setminus\mit\Gamma}^{}\frac{M_p}{2} \left[\alpha_\mrm{m} \mathrm{Tr}(\mbf{\eps_e}) - \zeta\right]^2 \, \mathd V
    ,
\end{equation}
where $\alpha_\mrm{m}$ and $M_p$ are Biot's coefficient and Biot's modulus, and $\zeta$ is the variation of fluid content~\cite{biot1962mechanics, mikelic2015quasi}.
Neglecting the effect of thermal expansion of the pore fluid, we have
\begin{equation}
    \label{eq:simi_incremental_content_of_fluid}
    \zeta =  \alpha_\mrm{m} \mathrm{Tr}(\mbfs{\eps_e}) + \frac{p}{M_p}.
\end{equation}

Assuming that the mechanical or hydraulic energy does not contribute to the thermal energy~\cite{suh2021asynchronous}, we have
\begin{equation}
    \label{eq:thermal_contribution}
    W_{\mathrm{thermo}}(T) := \int_{\mit{\Omega}\setminus\mit\Gamma}^{}\psi_T(T) \, \mathd V = \int_{\mit{\Omega}\setminus\mit\Gamma}^{}(\rho c)_m\left[(T - T_{\mathrm{ref}}) - T \ln\left(\frac{T}{T_{\mathrm{ref}}}\right)\right] \, \mathd V.
\end{equation}
With the local thermal equilibrium assumption, $T$ represents the unified temperature in the domain, and $T_{\mathrm{ref}}$ is the reference temperature. 
$(\rho c)_m$ is the equivalent heat storage coefficient, and for porous medium, it is given by
\begin{equation}
    \label{eq:heat_storage}
    (\rho c)_m = \phi c_f\rho_{p,f} + (1 - \phi) c_s\rho_{p,s}
    ,
\end{equation}
where $c$ is the specific heat, $\rho$ is the density, and $\phi$ is the porosity.

Substituting Eq.~\eqref{eq:thermo_elastic_strain_energy},~\eqref{eq:pore_fluid_energy},~\eqref{eq:thermo_elastic_strain_energy_density},~\eqref{eq:thermal_contribution} into the elastic strain energy in Eq.~\eqref{eq:thermo_poroelastic_strain_energy_total} yields
\begin{equation}
    \begin{aligned}
        \label{eq:thermo_poroelastic_strain_energy}
        W(\mbf{u}, p, T) 
        & = \int_{\mit{\Omega}\setminus\mit\Gamma}^{}\psi_e(\mbf{u}, T) \, \mathd V + \int_{\mit{\Omega}\setminus\mit\Gamma}^{}\frac{M_p}{2} \left[\alpha \mathrm{Tr}(\mbfs{\eps_e}) - \zeta\right]^2 \, \mathd V\\
        & + \int_{\mit{\Omega}\setminus\mit\Gamma}^{}(\rho c)_m\left[(T - T_{\mathrm{ref}}) - T\ln\left(\frac{T}{T_{\mathrm{ref}}}\right)\right] \, \mathd V\\
        & = \int_{\mit{\Omega}\setminus\mit\Gamma}^{}\frac{1}{2} \mathbb{C}_{m} : \mbfs{\eps}_e(\mbf{u}) : \mbfs{\eps}_e(\mbf{u}) \, \mathd V + \int_{\mit{\Omega}\setminus\mit\Gamma}^{}\frac{M_p}{2} \left[\alpha \mathrm{Tr}(\mbf{\eps_e}) - \zeta\right]^2 \, \mathd V\\
        & + \int_{\mit{\Omega}\setminus\mit\Gamma}^{}(\rho c)_m\left[(T - T_{\mathrm{ref}}) - T\ln\left(\frac{T}{T_{\mathrm{ref}}}\right)\right] \, \mathd V
        .
    \end{aligned}
\end{equation}

From Eq.~\eqref{eq:thermo_poroelastic_strain_energy}, we obtain the total Cauchy stress tensor $\mbf{\sigma}$, which can be decomposed into effective stress tensor $\sigma_{\mathrm{\mathrm{eff}}}$, pore pressure and thermal stress\footnote{We follow the engineering sign convention for stresses (positive for tension).}:
\begin{equation}
    \label{eq:thermo_poroelastic_stress}
    \begin{aligned}
    \frac{\partial W}{\partial\eps} 
    & :=  \boldsymbol{\sigma} \\
    & = \boldsymbol{\sigma}_{\mathrm{eff}} - \alpha_\mrm{m} p\mbf{I} - 3\beta K\Delta T\mbf{I} \\
    & = \mathbb{C}_{m}:\boldsymbol{\eps}(\mbf{u}) - \alpha_\mrm{m} p\mbf{I} - 3\beta K\Delta T\mbf{I}
    \end{aligned}
    ,
\end{equation}
where $K$ is the drained bulk modulus. 

Substituting Eq.~\eqref{eq:thermo_poroelastic_strain_energy} into Eq.~\eqref{eq:thermo_energy_functional} yields
\begin{equation}
    \begin{aligned}
        \label{eq:potential_energy}
        \mathcal{F}_\Gamma(\mbf{u}, \mit{\Gamma}, p, T) 
        & = W(\mbf{u}, p, T) + \int_{\mit{\Gamma}}^{}G_c \, \mathd S - \int_{\mit{\Omega}\setminus\mit\Gamma}^{}\mbf{b}\cdot\mbf{u} \, \mathd V - \int_{\mit{\mathcal{C}}_\mrm{N}}^{}\mbf{\bar{t}}\cdot\mbf{u} \, \mathd V \\
        & = \int_{\mit{\Omega}\setminus\mit\Gamma}^{}\frac{1}{2} \mathbb{C}_{m} : \mbfs{\eps}_e(\mbf{u}) : \mbfs{\eps}_e(\mbf{u}) \, \mathd V + \int_{\mit{\Omega}\setminus\mit\Gamma}^{}\frac{M_p}{2} \left[\alpha \Tr{\mbfs{\eps_e}} - \zeta\right]^2  \mathd V \\
        & +\int_{\mit{\Omega}\setminus\mit\Gamma}^{}(\rho c)_m\left[(T - T_{\mathrm{ref}}) - T \ln\left(\frac{T}{T_{\mathrm{ref}}}\right)\right] \, \mathd V + \int_{\mit{\Gamma}}^{}G_c \, \mathd S - \int_{\mit{\Omega} \setminus \mit{\Gamma}}^{}\mbf{b}\cdot\mbf{u} \, \mathd V - \int_{\mit{\mathcal{C}}_\mrm{N}}^{}\mbf{\bar{t}}\cdot\mbf{u} \, \mathd V
    \end{aligned}
    ,
\end{equation}
where $\mbf{\bar{t}}$ contains the thermal stress as
\begin{equation}
    \label{eq:boundary_stress}
    \mbf{\bar{t}} = \Big(\mathbb{C}_{m}:\mbfs{\eps}(\mbf{u}) - 3\beta K\Delta T\mbf{I}\Big)\cdot\mbf{n}
\end{equation}
with $\mbf{n}$ being the outward unit normal vector to the outer surface $\mit{\mathcal{C}_\mrm{N}}$. 
To sum up, the potential energy consists of the following parts: mechanical term, i.e., the strain energy from the linear elastic constitutive relationship; hydrothermal work in the domain; surface energy of fracture; the work of body force and the external load from the traction force.

\subsection{Mechanical equilibrium with phase-field approximation}
To alleviate the implementational difficulties of the sharp fractures, we apply the now-classical phase-field approach~\cite{bourdin2000numerical,bourdin2008variational}. 
Bourdin et al. (2000)~\cite{bourdin2000numerical} introduced a phase-field variable $\upsilon$ that represents a state of the material from intact $(\upsilon = 1)$ to fully broken $(\upsilon = 0)$ and regularized the Francfort-Marigo energy functional as
\begin{equation}
    \label{eq:pf_energy_functional}
    \mathcal{E}(\mbf{u}, \upsilon) := \int_{\mit{\Omega}}^{}\psi(\mbf{u}, \upsilon) \, \mathd V + \int_{\mit{\Omega}}^{}\frac{G_c}{4c_n}\left[\frac{(1-\upsilon)^n}{\ell} + \ell\nabla\upsilon\cdot\nabla\upsilon\right]\, \mathd V - \int_{\mit{\Omega}\setminus\mit\Gamma}^{}\mbf{b}\cdot\mbf{u} \, \mathd V - \int_{\mit{\mathcal{C}}_\mrm{N}}^{}\mbf{\bar{t}}\cdot\mbf{u} \, \mathd V
    ,
\end{equation}
where $c_n$ is the normalizing parameter defined as $c_n = \int_{0}^{1}(1 - s)^{n/2} \, \mathd S$~\cite{mesgarnejad2015validation,tanne2018crack}. 
For $n = 1$, the model is called \ATone{} model, and for $n = 2$, \ATtwo{} model~\cite{pham2011gradient}. 
Also, $\ell$ is the characteristic parameter with the dimension of a length that controls the phase-field profile transition. 
The strain energy density $\psi(\mbf{u}, \upsilon)$ in Eq.~\eqref{eq:pf_energy_functional} acknowledges the phase-field variable $\upsilon$ and is continuous over $\mit{\Omega}$ for integration.

Now, extending the regularized functional in Eq.~\eqref{eq:pf_energy_functional} to the thermo-poro-elastic medium, we can write
\begin{equation}
    \begin{aligned}
        \label{eq:pf_thermo_poroelastic_energy_functional}
        &\mathcal{F}(\mbf{u}, p, T, \upsilon) := \int_{\mit{\Omega}}^{}\psi_e(\mbf{u}, T, \upsilon) \, \mathd V + \int_{\mit{\Omega}}^{}\psi_f(\mbf{u}, \zeta, \upsilon) \, \mathd V + \int_{\mit{\Omega}}^{}\psi_T(T, \upsilon) \, \mathd V\\
        & + \int_{\mit{\Omega}}^{}\frac{G_c}{4c_n}\left[\frac{(1-\upsilon)^n}{\ell} + \ell\nabla\upsilon\cdot\nabla\upsilon\right]\, \mathd V - \int_{\mit{\Omega}}^{}\mbf{b}\cdot\mbf{u} \, \mathd V - \int_{\mit{\mathcal{C}}_\mrm{N}}^{}\mbf{\bar{t}}\cdot\mbf{u} \, \mathd V.
    \end{aligned}
\end{equation}
The total stress can be found by taking the derivative of $\mathcal{F}$ with respect to $\mbfs{\eps}_e$ as
\begin{equation}
    \label{eq:pf_thermo_poroelastic_stress}
    \sigma = \mathbb{C}_{\mathrm{eff}}(\upsilon):\mbf{\eps}(\mbf{u}) - \alpha(\upsilon) p\mbf{I} - 3\beta K\Delta T\mbf{I}
    ,
\end{equation}
where $\mathbb{C}_{\mathrm{eff}}(\upsilon)$ and $\alpha(\upsilon)$ are the effective stiffness tensor and effective Biot's coefficient.
Comparing with Eq.~\eqref{eq:thermo_poroelastic_stress}, we note that $\mathbb{C}_{\mathrm{eff}}(\upsilon)$ and $\alpha(\upsilon)$ represent both intact and the fractured medium while $\mathbb{C}_\mrm{m}$ and $\alpha_\mrm{m}$ represent only the intact material. 

The expressions of $\mathbb{C}_\mrm{eff} (v)$ and $\alpha(\upsilon)$ can be derived from the strain energy decomposition scheme~\cite{YOU2023116305}.
The elastic strain energy density is generally split into the positive (tension) and negative (compression) parts, $\psi = \psi_+ + \psi_-$, to differentiate the contributions between positive (tension) and negative (compression) parts to the damage evolution. 
With a degradation function $g(\upsilon)$\footnote{In this study, we employed $g(\upsilon) = (1 - k)\upsilon^2 + k$ where $k$ is a phase-field parameter representing residual stiffness, which keeps the system of equations well-conditioned for the partly-broken state~\cite{martinez2018phase}.}, $\psi_e (\mbf{u},T, \upsilon)$ is rewritten as~\cite{amor2009regularized,freddi2010regularized,miehe2010phase}
\begin{equation}
    \label{eq:pf_strain_energy}
    \begin{aligned}
    \psi_e (\mbf{u},T,\upsilon) 
    &= g(\upsilon) \psi_+(\mbf{u}) + \psi_-(\mbf{u}) \\
    &= \frac12 g(\upsilon) \mathbb{C}_+:\mbfs{\eps}_e : \mbfs{\eps}_e + \frac12 \mathbb{C}_-:\mbfs{\eps}_e : \mbfs{\eps}_e \\
    &= \frac12 \mathbb{C}_{\mathrm{eff}}(\upsilon):\mbfs{\eps}_e : \mbfs{\eps}_e
    ,
    \end{aligned}
\end{equation}
where $\mathbb{C}_{\mathrm{eff}}(\upsilon) = g(\upsilon) \mathbb{C}_+ + \mathbb{C}_-$.
Similarly, Biot's coefficient evolves with the phase-field as the drained bulk modulus is degraded with damage as~\cite{YOU2023116305}
\begin{equation}
    \label{eq:pf_Biot_coeff}
    \alpha(\upsilon) = 1 - \frac{K_{\mathrm{eff}}(\upsilon)}{K_s}.
\end{equation}

The expressions of $\mathbb{C}_{\mathrm{eff}}(\upsilon)$ and $K_{\mathrm{eff}}(\upsilon)$ depend on the specific energy decomposition model applied as discussed in~\cite{YOU2023116305}.
We employed the volumetric-deviatoric (V-D) energy split~\cite{amor2009regularized} in this study and then the strain energy is split into~\cite{amor2009regularized}
\begin{equation}
    \label{eq:V_D_method}
    \psi^{\mathrm{vd}}(u,v)
    =g(v)\left[\frac{K_\mrm{m}}{2}\langle\mathrm{Tr}(\mbfs{\eps}_e)\rangle_{+}^{2}+\mu\mbfs{\eps}_e:\mathbb{K}:\mbfs{\eps}_e\right]+\frac{K_\mrm{m}}{2}\langle\mathrm{Tr}(\mbfs{\eps}_e)\rangle_{-}^{2}
    ,
\end{equation}
where $\langle\cdot\rangle$ denotes the Macaulay brackets defined as $\langle\cdot\rangle_{\pm} = (|\cdot|\pm\cdot)/2$. 
Accordingly, we can derive the tangential stiffness ($\partial^2 \psi^\mrm{vd} / \partial \boldsymbol{\eps}^2$) as
\begin{equation}
    \label{eq:C_eff}
    \begin{aligned}
    \mathbb{C}_{\mathrm{eff}} (v)
   & = 3\left[g(\upsilon)H(\mathrm{Tr}(\mbfs{\eps}_e))+H(\mathrm{Tr}(-\mbfs{\eps}_e))\right]K_\mrm{m}\mathbb{J}+2g(\upsilon)\mu\mathbb{K} \\
   & = 3K_{\mathrm{eff}}\mathbb{J} + 2g(\upsilon)\mu\mathbb{K} 
    \end{aligned}
    ,
\end{equation}
with
\begin{equation}
    \label{eq:K_eff}
    K_{\mathrm{eff}} = \left[g(\upsilon)H(\mathrm{Tr}(\mbfs{\eps}_e))+H(\mathrm{Tr}(-\mbfs{\eps}_e))\right]K_\mrm{m}
    ,
\end{equation}
where $H(\cdot)$ is the Heaviside step function defined as
\begin{equation}
    \label{eq:Heaviside}
    H(x):=
    \begin{cases}
    0, \quad &x<0\\
    1,\quad &x\geq0
    \end{cases}
    .
\end{equation}

Substituting Eq.~\eqref{eq:K_eff} into Eq.~\eqref{eq:pf_Biot_coeff}, degraded Biot's coefficient based on its initial value $\alpha_\mrm{m}$ is
\begin{equation}
    \label{eq:pf_Biot_coeff_final}
    \begin{aligned}
    \alpha(v)& =1-\left[g(\upsilon)H(\mathrm{Tr}\left(\mbfs{\eps}_e\right))+H(\mathrm{Tr}\left(-\mbfs{\eps}_e\right))\right]\frac{K_\mrm{m}}{K_s}  \\
    &=1-\left[g(v)H(\mathrm{Tr}\left(\mbfs{\eps}_e\right))+H(\mathrm{Tr}\left(-\mbfs{\eps}_e\right))\right](1-\alpha_\mrm{m}).
    \end{aligned}
\end{equation}
Eq.~\eqref{eq:pf_Biot_coeff_final} shows that Biot's coefficient depends both on the damage and the crack opening. When the fracture opens, i.e., $\mathrm{Tr}\left(\mbfs{\eps}\right)\geq 0$, it is enhanced to $\alpha = 1 - g(\upsilon) + g(\upsilon)\alpha_\mrm{m} = 1$. 
Otherwise, the fracture is closed, and $\alpha_m\mrm{m}= 1 - (1 - \alpha_\mrm{m}) = \alpha_\mrm{m}$.

\subsection{Fluid flow model}
Considering the equivalent properties over $\mit \Omega$, the mass balance is given as
\begin{equation}
    \frac{\partial}{\partial t}\left(\alpha (\upsilon) \nabla\cdot\mbf{u} + \frac{p}{M_p (\upsilon)} - \frac{T}{M_T (\upsilon)}\right) +\nabla\cdot\left(\mbf{q}_f\right)=Q_{f} \quad \mrm{in} \quad \Omega
    ,
    \label{eq:mass_balance_equation}
\end{equation}
where $M_T(\upsilon)$ is the effective thermal storage coefficient that describes the thermal expansion in the incremental content of pore fluid. 
Before we derive the expressions of $M_p(\upsilon)$ and $M_T(\upsilon)$, we first discuss how the permeability tensor, $\mbf{K}$, is obtained and the porosity $\phi$ is updated in our model. 
The fluid flux $\mbf{q}_f$ is given by Darcy's law as
\begin{equation}
    \mbf{q}_f=-\frac{\mbf{K}}{\mu}\nabla(p + \gamma_f z)
    \quad \mrm{in} \quad \Omega
    ,
    \label{eq:fluid_flux}
\end{equation}
with $\mu$ the fluid viscosity, $\gamma_f$ the specific weight of water and $z$ the vertical coordinate, and $Q_f$ the source term. 
In the following derivation, the gravity can be neglected without loss of generality. 

For permeability enhancement by fractures, we apply a formulation of anisotropic permeability which implicitly takes into account for the Poiseuille-type flow in fractures~\cite{miehe2016phase,miehe2015minimization}
\begin{equation}
    \mbf{K}=\mbf{K_m}\mbf{I}+(1 - \upsilon)^\xi\frac{w^2}{12}\left(\mbf{I}-\mbf{n}_{\mit{\Gamma}}\otimes\mbf{n}_{\mit{\Gamma}}\right) 
    \label{eq:enhanced_permeability}
    ,
\end{equation}
where $\mbf{K_m}$ is the isotropic permeability, $\xi\geq 1$ is a weighting exponent, and the fracture has an aperture of $w$ with the normal vector $\mbf{n}_{\mit{\Gamma}}$ along the interface. 
The weighting exponent $\xi$ controls the intensity of permeability enhancement to ensure numerical stability. 

One may approximate the crack normal vector ($\mbf{n}_{\mit{\Gamma}}$) from the gradient of the phase-field $\nabla\upsilon$, but $\nabla\upsilon$ deviates from the normal direction near the crack tip~\cite{bourdin2012variational}. 
Moreover, $\nabla\upsilon$ is not defined on a fully broken element where $\upsilon = 0$. 
To avoid these problems, we decompose a strain tensor into principal strains as
\begin{equation}
    \bm{\eps}=\sum_{i=1}^{3}\eps_i \mbf e_i 
    \label{eq:decomposed_strain_tensor}
\end{equation}
where $\eps_i$ are the principal strains ($\eps_1 > \eps_2 > \eps_3$) and $\mbf e_i$ are the associated eigenvectors. 
Considering an internally pressurized crack, the crack normal deformation must be dominant, and thus we assign
\begin{equation}
    \mbf{n}_{\mit{\Gamma}} = \mbf e_1
    .
    \label{eq:normal_vector}
\end{equation}
This approach enables simple code parallelization because the computation is performed at each integration point, while some other methods like line integral~\cite{bourdin2012variational} or level-set~\cite{lee2017iterative} based approach require a global phase-field profile.

To estimate the fracture width ($\omega$), Miehe et al. (2015)~\cite{miehe2015minimization} proposed:
\begin{equation}
    \label{eq:frac_width0}
    \omega = h_e\eps_{\mathrm{vol}}
    ,
\end{equation}
where $h_e$ is the corresponding element size, and $\eps_{\mathrm{vol}}$ represents the volumetric strain. 
Eq.~\eqref{eq:frac_width0} can approximate the fracture width accurately as reported in~\cite{YOU2023116305} because the volumetric strain is almost identical to the strain in the crack normal opening direction for fluid-driven fractures.
However, this approximation can lead to spurious crack opening profiles when thermal stress is present. 
Instead, in this study, we approximate the fracture width using the maximum principal strain:
\begin{equation}
    \label{eq:frac_width}    
    \omega = h_e\eps_1.
\end{equation}
For further discussions, we refer to \ref{C}, where we compared crack opening profiles from Eq.~\eqref{eq:frac_width0} with those from Eq.~\eqref{eq:frac_width}.

Regarding the porosity update, existing studies undertook various ways to account for evolving phase-field (damage). 
Lee and Wick (2017)~\cite{lee2017iterative} first proposed a linear indicator function with threshold values to define the effective porosity distinguishing between the matrix and fracture
\begin{equation}
    \label{eq:pf_porosity_C}
    \phi(\upsilon) = \mathcal{X}_R\left(\frac{1}{M}p_R + \alpha\nabla\cdot\mbf{u}\right) + \mathcal{X}_F\left(c_fp_f\right)
\end{equation}
where $\mathcal{X}_R \in [0,1]$ and $\mathcal{X}_F\in [0,1]$ are two linear indicator functions for the reservoir domain and the fracture domain. 
$\mathcal{X}_F(\upsilon)$ is 0 and $\mathcal{X}_R(\upsilon)$ is 1 in the reservoir domain while $\mathcal{X}_F(\upsilon)$ is 1 $\mathcal{X}_R(\upsilon)$ is 0 in the fracture domain.
And $\mathcal{X}_F(\upsilon)$ and $\mathcal{X}_R(\upsilon)$ vary linearly in the transition zone.
This concept has been used in may works~\cite{guo2024reactive, xing2023hydro, jammoul2022phase, feng2016numerical, zhou2018phase}. 
Similarly, Li et al. (2021)~\cite{li2021phase} linearly interpolated the matrix porosity ($\phi_m$) and the fracture porosity (1) as
\begin{equation}
    \label{eq:pf_porosity_D}
    \phi(\upsilon) = \mathcal{X}_R\phi_m + \mathcal{X}_F
    .
\end{equation}
Zhou et al.~\cite{zhou2018phase, zhou2019phase} also applied a similar expression without considering the fluid compressibility
\begin{equation}
    \label{eq:pf_porosity_E}
    \phi(\upsilon) = \mathcal{X}_R\phi_m.
\end{equation}
In \cite{yi2024coupled}, the initial porosity $\phi_m$ is degraded with the phase field $\upsilon$ as
\begin{equation}
    \label{eq:pf_porosity_A}
    \phi(\upsilon) = \phi_m + (1 - \phi_m)(1 - (1 - \upsilon)^2)
    .
\end{equation}
Suh and Sun~\cite{suh2021asynchronous} additionally considered the effect of volumetric strain:
\begin{equation}
    \label{eq:pf_porosity_B}
    \phi(\upsilon) = \phi_m + g(\upsilon)(1 - \phi_m)(1 - \nabla\cdot\mbf{u}).
\end{equation}
Lastly, You and Yoshioka (2023)~\cite{YOU2023116305} updated the porosity in a similar manner to Biot's coefficient (Eq.\eqref{eq:pf_Biot_coeff_final}) as
\begin{equation}
    \label{eq:pf_porosity_0}
    \phi_0(\upsilon) = 1-\left[g(\upsilon)H(\mathrm{Tr}\left(\mbfs{\eps}_e\right))+H(\mathrm{Tr}\left(-\mbfs{\eps}_e\right))\right](1-\phi_{m}).
\end{equation}

All of these approaches consider the transition of porosity from the matrix ($\phi_\mrm{m}$) to fracture ($\phi = 1$) relying on the phase-field profile, $\upsilon (\mbf{x})$. 
While the phase-field profile, $\upsilon (\mbf{x})$, is a diffused representation of crack, it does not physically represent the porosity at the corresponding location $\mbf{x}$.
This is obvious because of the impacts of the characteristic length $\ell$ on the phase-field transition.
The longer $\ell$, the wider the transition zone.
However, the porosity profile should not change with $\ell$.

\begin{figure}[H]
		\centering
		\includegraphics[scale=0.4]{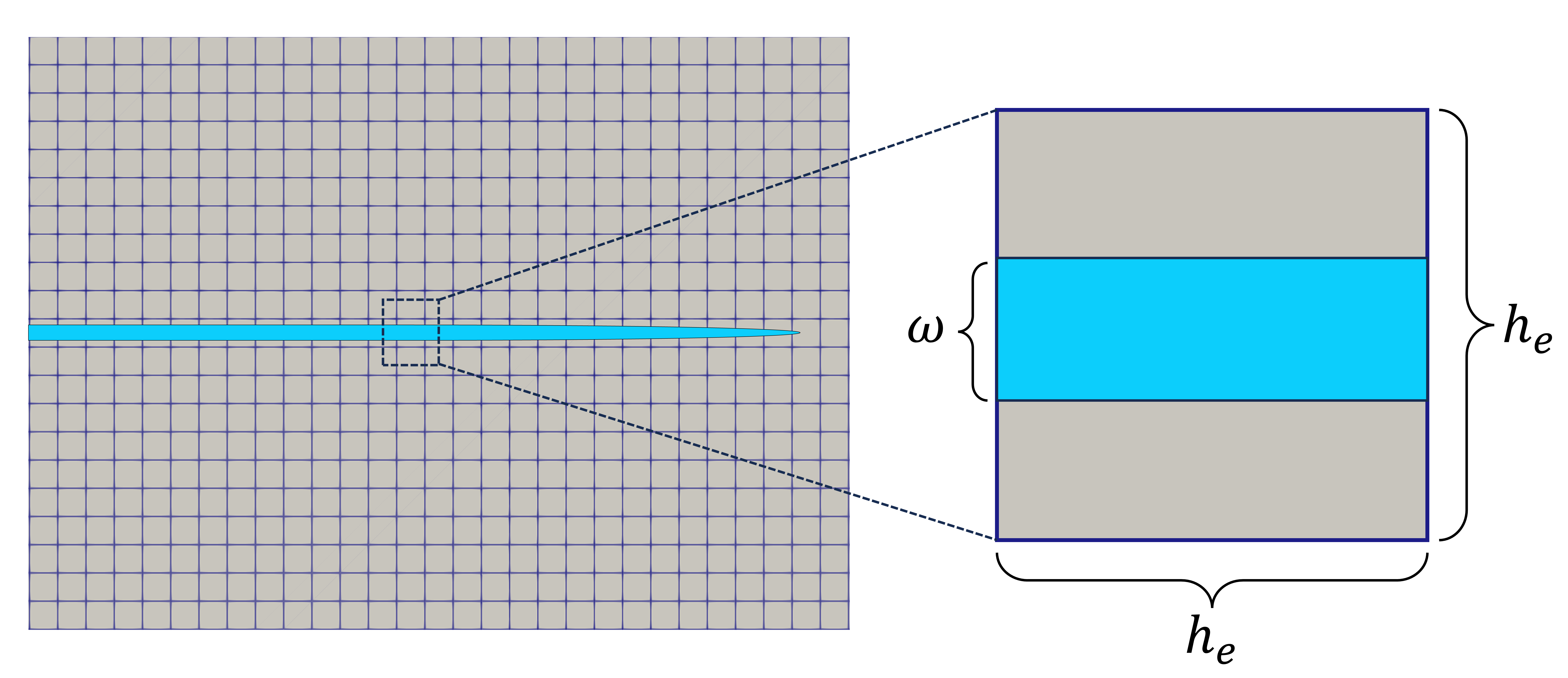}
		\caption{\centering{Physical crack opening within quadrilateral elements.}}
		\label{fig:pf_porosity}
\end{figure}

Considering that a hydraulic fracture width is generally in the order of 10$^{-2}$~m, it is likely smaller than a computational element size $h_e$, which is typically in the order of 10$^{-1}\sim $10$^{0}$~m.
Therefore, one can consider that a physical crack opening is well contained within one element as depicted in Fig.~\ref{fig:pf_porosity}. 
Accordingly, we can write the porosity change in a 2D quadrilateral element with the edge length $h_e $ containing a line crack parallel to the edge as
\begin{equation}
    \label{eq:pf_porosity}
    \phi_1 = \phi_m + \frac{\omega h_e}{{h_e}^2} = \phi_m + \frac{\omega}{h_e}.
\end{equation}
Recalling Eq.~\eqref{eq:frac_width}, the porosity is given as:
\begin{equation}
    \label{eq:pf_porosity_final}
    \phi_1(\eps) = \phi_m + \eps_1.
\end{equation}

Using this porosity expression and effective Biot's coefficient, we can write effective Biot's modulus ($M_p (\upsilon)$) and the average thermal expansion coefficient ($M_T (\upsilon)$) as
\begin{equation}
    \frac{1}{M_p(\upsilon,\eps)}=\phi(\eps) c_{f}+\frac{\alpha(\upsilon)-\phi(\eps)}{K_{s}},
    \label{eq:pf_storage_term}
\end{equation}
\begin{equation}
    \frac{1}{M_T(\upsilon,\eps)}=\phi(\eps) \alpha_f + 3\beta\left(\alpha(\upsilon)-\phi(\eps)\right).
    \label{eq:pf_average_thermal_expansion}
\end{equation}
where $c_f$ is the fluid compressibility, $K_s$ is the solid phase's intrinsic bulk modulus, and $\alpha_f$ is the fluid volumetric thermal expansion coefficient. 

\subsection{Heat transfer model}
Assuming the local thermal equilibrium and neglecting thermal effects due to deformation, we obtain the energy conservation equation in the porous medium as
\begin{equation}
         (\rho c)_m\frac{\partial T}{\partial t} + \nabla\cdot(\mbf{q}_T) =Q_T
         ,
    \label{eq:energy_conservation_equation}
\end{equation}
where $Q_T$ is the thermal source term.
The total heat flux $\mbf{q}_T$ in Eq.~\eqref{eq:heat_flux} is decomposed into the advective and the conductive terms using Fourier's law:
\begin{equation}
         \mbf{q}_T = \rho_f \mbf{q}_f c_f T - \mbf{\lambda}_{\mathrm{eff}}\nabla T
         ,
    \label{eq:heat_flux}
\end{equation}
with the effective thermal conductivity given as
\begin{equation}
         \mbf{\lambda}_{\mathrm{eff}}(\eps)=\phi(\eps)\mbf{\lambda}_f + (1-\phi(\eps))\mbf{\lambda}_s.
    \label{eq:pf_effective_thermal_conductivity}
\end{equation}
Similarly, the effective heat storage coefficient (Eq.~\eqref{eq:heat_storage}) is rewritten as
\begin{equation}
    \label{eq:pf_heat_storage}
    (\rho c)_m(\eps) = \phi(\eps) c_{p,f}\rho_f + (1 - \phi(\eps)) c_{p,s}\rho_s
    .
\end{equation}
Substituting Eq.\eqref{eq:heat_flux} into Eq.~\eqref{eq:energy_conservation_equation}, we have
\begin{equation}
    \begin{aligned}
        &(\rho c)_{m}\frac{\partial T}{\partial t} + c_{p,f}\rho_f\mbf{q}_f\cdot\nabla T - \nabla\cdot\mbf{\lambda}_{\mathrm{eff}}\nabla T = Q_T.
        \label{eq：new_heat_transfer_model}
    \end{aligned}
\end{equation}

%% file: numeric.tex
In a quasi-static setting, we obtain a solution pair of $(\mbf{u}_i, \upsilon_i)$ that minimizes $\mathcal{F}$ (Eq.~\eqref{eq:pf_thermo_poroelastic_energy_functional}) each time as 
\begin{equation}
	\label{eq:glob_min}
	(\mathbf{u}_i, \upsilon_i) =
	\argmin \left\lbrace  \mathcal{F}(\mathbf{u},\upsilon, p, T) : \mathbf{u} \in \mathcal{U}(t_i), \upsilon \in \mathcal{\upsilon}(t_i, \upsilon_{i-1}) \right\rbrace, 
\end{equation}
where $\mathcal{U}$ is the kinematically admissible displacement set:
\begin{equation}
	\label{eq:_adm}
	\mathcal{U}(t_i) = \left\{\mathbf{u} \in H^1 (\Omega) : \mathbf{u} = 0 \quad \text{on} \quad \partial_N \Omega \right\}.
\end{equation}
The kinematically admissible set of $\upsilon$ requires an irreversible condition as:
\begin{equation}
\label{eq:v_adm}
\mathcal{V}(t_i,\upsilon_{i-1}) = \left\{\upsilon \in H^1 (\Omega) : 0 \le \upsilon(x) \le \eta \;\; \forall x \right\}.
\end{equation}
where
$$ \eta =
\left\{
\begin{array}{ll}
		1 &\quad \text{if} \quad  \upsilon_{i-1}(x) \ge \upsilon_{ir} \\
		\upsilon_{i-1}(x) &\quad \text{otherwise}
\end{array}
	\right .
$$
and $\upsilon_{ir}$ is the irreversible threshold $\in [0,1]$ ({e.g.}~0.05).
Following~\cite{bourdin2000numerical}, we apply an alternate minimization scheme which minimizes $\mathcal{F}$ with respect to $\mbf{u}$ while fixing $\upsilon$ and then minimizes $\mathcal{F}$ with respect to $\upsilon$ while fixing $\mbf{u}$ with the irreversible condition. 
Therefore, the variations of $\mathcal{F}$ with respect to $\mbf{u}$ and $\upsilon$ are given:
\begin{equation}
    \begin{aligned}
    	\delta_\mbf{u}\mathcal{F} = &
        \int_{\mit\Omega}g(\upsilon)\mbf{\sigma}_e:\delta\mbf{\eps} \, \mathd V - \int_{\mit\Omega}\alpha(\upsilon)p\delta\eps \, \mathd V-\int_{\mathcal{C}_N}\bar{\mbf{t}}\cdot\delta \mbf{u}\, \mathd S - \int_{\mit{\Omega}}^{}\mbf{b}\cdot\mbf{u} \, \mathd V\\
    	 = &
    	-\int_{\mit\Omega}\nabla\cdot(g(\upsilon)\mbf{\sigma}_e)\cdot\delta\mbf{u} \, \mathd V +\int_{\partial_{N}\mit\Omega}g(\upsilon)\mathbb{C}_{eff}:\mbf{\eps}_e\cdot\mbf{n}\cdot\delta\mbf{u}\, \mathd S -\int_{\mathcal{C}_N}\alpha(\upsilon)p\mbf{I}\cdot\mbf{n}\cdot\delta\mbf{u}\, \mathd S\\
    	&+\int_{\mit\Omega}\nabla\cdot(\alpha(\upsilon)p\mbf{I})\cdot\delta\mbf{u}\, \mathd V - \int_{\mathcal{C}_N}\bar{\mbf{t}}\cdot\delta \mbf{u}\, \mathd S - \int_{\mit{\Omega}}^{}\mbf{b}\cdot\delta\mbf{u} \, \mathd V\\
         =  &\mbf{0},
    	\label{eq:u_varia}
    \end{aligned}
\end{equation}
and
\begin{equation}
    \begin{aligned}
    	\delta_\upsilon\mathcal{F} =&
    	\int_{\mit\Omega}2(1-k)\upsilon\psi_+(\mbf{u})\delta\upsilon \, \mathd V +\int_{\mit\Omega}\frac{p^2}{2}\frac{\partial 1/M_p(\upsilon,\mbf{\eps})}{\partial\upsilon}\delta\upsilon \, \mathd V\\
        &+ \frac{G_{c}}{4c_{n}}\int_{\mit\Omega}[-\frac{n}{\ell}(1-\upsilon)^{n-1}\delta\upsilon+2\ell\nabla\upsilon\cdot\nabla\delta\upsilon]\, \mathd V \\
    	=&
        \int_{\mit\Omega}2(1-k)\upsilon\psi_+(\mbf{u})\delta\upsilon dV + \int_{\mit\Omega}\frac{p^2}{2}\frac{\partial 1/M_p(\upsilon,\mbf{\eps})}{\partial\upsilon}\delta\upsilon \, \mathd V\\
        & + \frac{G_{c}}{4c_{n}}\int_{\mit\Omega}[-\frac{n}{\ell}(1-\upsilon)^{n-1}\delta\upsilon-2\ell\Delta\upsilon\delta\upsilon]\, \mathd V + \frac{G_{c}}{4c_{n}}\int_{\partial\mit\Omega}2\ell\nabla\upsilon\cdot\bm{n}\delta\upsilon \, \mathd S \\
        & = \mbf{0}
        \label{eq:pf_varia}
        .
    \end{aligned}
\end{equation}
Note that the thermal energy $W_{thermo}$ cancels because the porosity does not change with the phase-field~\cite{suh2021asynchronous,dittmann2020phase,na2018computational}. 
The derivative of the reciprocal of Biot's modulus $M_p(\upsilon)$ in Eq.~\eqref{eq:pf_varia} with respect to $\upsilon$ writes
\begin{equation}
    \label{de_pf_Biot_modu}
    \frac{1/M_p(\upsilon,\mbf{\eps})}{\partial\upsilon} = \frac{2\eps_\mathrm{vol}}{p}\upsilon(1-k)H( \Tr{\eps})\left(1 - \alpha_{\mathrm{m}}\right)
    ,
\end{equation}
and its derivation details are given in \ref{A}. 

From Eqs.~\eqref{eq:u_varia} and~\eqref{eq:pf_varia}, we arrive at the strong forms of the coupled problem.
For the mechanical deformation, we have
\begin{align}
\begin{dcases}
\nabla\cdot[g(\upsilon)\mathbb{C}_+:\mbf{\eps}_e -\alpha(\upsilon)p\bm{I}] + b = 0 & \text{in } \quad \mit{\Omega} \\
g(\upsilon)\mathbb{C}_+:\mbf{\eps}_e\cdot\mbf{n}-\bar{\mbf{t}} = 0 \quad &\text{on} \quad \mathcal{C}_{N} \\        g(\upsilon)\mathbb{C}_+:\mbf{\eps}_e\cdot\mbf{n} = 0 \quad &\text{on} \quad \partial\mit\Gamma \\
\end{dcases}
\label{eq:disfield}
,
\end{align}
and for the phase-field evolution,
\begin{align}
\begin{dcases}
(1-k)\upsilon\mathbb{C}_+:\mbf{\eps}_e:\mbf{\eps}_e+\frac{G_{c}}{4c_{n}}\left[-\frac{n}{\ell}(1-\upsilon)^{n-1}-2\ell\Delta\upsilon\right] +\frac{p^2}{2}\frac{\partial 1/M_p(\upsilon)}{\partial\upsilon} = 0 & \text{in } \quad \mit{\Omega} \\
\nabla\upsilon\cdot\bm{n}=0 & \text{on} \quad \partial_N\mit\Omega
\end{dcases}
\label{eq:pffield}
.
\end{align}

Multiplying Eqs.~\eqref{eq:disfield} and \eqref{eq:pffield} with the weighting functions $\mbf{w}_u\in H^1$ and $\mbf{w}_\upsilon\in H^1$ and integrating over $\mit{\Omega}$, we obtain the weak forms:
\begin{equation}
    \label{weak_form_dis}
    \int_{\mit\Omega}\nabla\mbf{w}_u\cdot\left[g(\upsilon)\mbf{\sigma}_e -\alpha(\upsilon)p\mbf{I}\right]\, \mathd V - \int_{\mit\Omega}\mbf{b}\cdot\mbf{w}_u\, \mathd V - \int_{\mathcal{C}_N}\bar{\mbf{t}}\cdot\mbf{w}_u\, \mathd S = \mbf{0}
\end{equation}
\begin{equation}
    \begin{aligned}
        \label{weak_form_pf}
        &\int_{\mit\Omega}\mbf{w}_\upsilon\left[2(1-k)\upsilon\psi_+(\mbf{u}) - \frac{p^2}{2}\frac{\partial 1/M_p(\upsilon)}{\partial\upsilon}\right]\, \mathd V - \int_{\mit\Omega}\mbf{w}_\upsilon\frac{G_{c}}{4c_{n}}\frac{n}{\ell}(1-\upsilon)^{n-1}\, \mathd V,\\
        &- \int_{\mathcal{C}_N}\frac{G_{c}}{2c_{n}}\ell\nabla \mbf{w}_\upsilon\cdot\nabla\upsilon\, \mathd V = \mbf{0}.
    \end{aligned}
\end{equation}
Similarly, with $\mbf{w}_p\in H^1$ and $\mbf{w}_T \in H^1(\Omega)$, the weak forms of the mass and energy balances are
\begin{equation}
    \begin{aligned}
        &\int_{\mit\Omega}^{}\frac{\partial}{\partial t}\left(\alpha(\upsilon)\nabla\cdot\mbf{u}+\frac{1}{M_p(\upsilon,\eps)}p-\frac{1}{M_T(\upsilon,\eps)}T\right)\mbf{w}_p \, \mathd V + \int_{\mit\Omega}^{}\frac{\mbf{K}}{\mu}\nabla p\cdot\nabla\mbf{w}_p \, \mathd V\\
        &= \int_{\mit\Omega} Q_f\mbf{w}_p \, \mathd V - \int_{\mathcal{C}_\mrm{N}}^{}q_n\mbf{w}_p \, \mathd S,
        \label{eq:pf_weak_form_of_fluid_flow_model}
    \end{aligned}
\end{equation}
and
\begin{equation}
    \begin{aligned}
        &\int_{\mit\Omega}^{}(\rho c)_m(\eps)\frac{\partial T}{\partial t}\mbf{w}_T \, \mathd V + \int_{\mit\Omega}^{}\mbf{w}_Tc_f\rho_f\mbf{q}_f\cdot\nabla T \, \mathd V +\int_{\mit\Omega}^{} \mbf{\lambda}_{\mathrm{eff}}(\eps)\nabla T\cdot\nabla\mbf{w}_T \, \mathd V\\
        & = \int_{\mit\Omega} Q_T\mbf{w}_p \, \mathd V - \int_{\mathcal{C}_\mrm{N}}^{}q_{Tn}\mbf{w}_T \, \mathd S.
        \label{eq:pf_weak_form_of_heat_transfer_model}
    \end{aligned}
\end{equation}

\subsection{The staggered solution scheme}
Eqs.~\eqref{weak_form_dis}, \eqref{weak_form_pf}, \eqref{eq:pf_weak_form_of_fluid_flow_model} and \eqref{eq:pf_weak_form_of_heat_transfer_model} form the non-linear system of partial differential equations describing thermo-hydro-mechanical phase-field coupled problem. 
We employed a staggered scheme based on~\cite{brun2020monolithic,wang2009parallel} to split the system into four sub-problems in a sequence of $\upsilon - (T-p-\mbf{u})$. 
Globally, we have a staggered loop between the process of $\upsilon$ and a paired process of $(T-p-\mbf{u})$.
Within the $(T-p-\mbf{u})$ loop, $T$, $p$, and $\mbf{u}$ are solved in a staggered manner until convergence.

Specifically, at time step $k$, the backward Euler scheme for the time derivative is
\begin{align}
\begin{dcases}
        &\frac{\partial(\nabla\cdot\mbf{u})}{\partial t} = \frac{\eps_\mathrm{vol}^{k} - \eps_\mathrm{vol}^{k-1}}{\Delta t} \\
		&\frac{\partial p}{\partial t} = \frac{p^k - p^{k-1}}{\Delta t} \\
        &\frac{\partial T}{\partial t} = \frac{T^k - T^{k-1}}{\Delta t} \\
\end{dcases}
.
\end{align}
Then, the $m^{th}\:(m\geq 1)$ iteration scheme reads
\vspace{20pt}

$\bullet$\:$\bm{\textbf{Step}\:1}$: Given\: $(T^{k,m-1}, p^{k,m-1}, \mbf{u}^{k,m-1})$, solve for\: $\upsilon^{k,m}$:
\begin{equation}
    \begin{aligned}
        \label{weak_form_pf_stagg}
        &\int_{\mit\Omega}\mbf{w}_\upsilon\left[2(1-k)\upsilon^{k,m}\psi_+(\mbf{u})^{k,m-1} - \frac{\left(p^{k,m-1}\right)^2}{2}\frac{\partial 1/M_p}{\partial\upsilon}\right]\, \mathd V\\
        &- \int_{\mit\Omega}\mbf{w}_\upsilon\frac{G_{c}}{4c_{n}}\frac{n}{\ell}(1-\upsilon^{k,m})^{n-1}\, \mathd V
        - \int_{\mathcal{C}_N}\frac{G_{c}}{2c_{n}}\ell\nabla \mbf{w}_\upsilon\cdot\nabla\upsilon^{k,m}\, \mathd V = \mbf{0}.
    \end{aligned}
\end{equation}
\vspace{20pt}

$\bullet$\:$\bm{\textbf{Step}\:2}$: Given\: $(p^{k,m-1}, \upsilon^{k,m})$, find\: $T^{k,m}$ independently:
\begin{equation}
    \begin{aligned}
        &\int_{\mit\Omega}^{}(\rho c)_m(\eps)\frac{T^{k,m}-T^{k-1}}{\Delta t}\mbf{w}_T \, \mathd V + \int_{\mit\Omega}^{}\mbf{w}_Tc_f\rho_f\mbf{q}_f^{k,m-1}\cdot\nabla T \, \mathd V +\int_{\mit\Omega}^{} \mbf{\lambda}_\mrm{eff}(\eps)\nabla T^{k,m}\cdot\nabla\mbf{w}_T \, \mathd V\\
        & = \int_{\mit\Omega} Q_T\mbf{w}_p \, \mathd V - \int_{\mathcal{C}_N}^{}q_{Tn}\mbf{w}_T \, \mathd S.
        \label{eq:pf_weak_form_of_heat_transfer_model_stagg}
    \end{aligned}
\end{equation}

    with 
\begin{equation}
    \mbf{q}_f^{k,m-1} = -\frac{\mbf{K}}{\mu}\nabla p^{k,m-1}
    \label{eq:velit}
\end{equation}
\vspace{20pt}

$\bullet$\:$\bm{\textbf{Step}\:3}$: Given\: $(T^{k,m}, \mbf{u}^{k,m-1}, p^{k,m-1}, \upsilon^{k,m})$, solve for\: $p^{k,m}$:
\begin{equation}
    \begin{aligned}
        &\int_{\mit\Omega} (\alpha \frac{\varepsilon_v(\mbf{u}^{k,m})-\varepsilon_v(\mbf{u}^{k-1})}{\Delta t}+\frac{1}{M_p}\frac{p^{k,m}-p^{k-1}}{\Delta t}-\frac{1}{M_T}\frac{T^{k,m}-T^{k-1}}{\Delta t})\psi_{p} dV + \int_{\mit\Omega}\frac{\bm{K}}{\mu}\nabla p^{k,m}\cdot\nabla\psi_{p} dV  \\
        & = \int_{\Omega} Q_{f}\psi_{p} dV-\int_{\partial_{N}\mit\Omega}q_{n}\psi_{p} dS.
        \label{eq:weak_form_of_fluid_flow_model_stagg}
    \end{aligned}
\end{equation}
We apply the fixed-stress splitting method proposed by~\cite{kim2011stability} for the stability of the fluid flow equation. 
Freezing the volumetric stress, we eliminate the volumetric strain in the last iteration from the mass balance:
\begin{equation}
\begin{aligned}
    &K_{\mathrm{eff}}\eps_\mathrm{vol}(\mbf{u}^{k,m})-\alpha p^{k,m}-3\alpha_sK_{\mathrm{eff}}(T^{k,m}-T_0)\\
    &=K_{\mathrm{eff}}\eps_\mathrm{vol}(\mbf{u}^{k,m-1})-\alpha p^{k,m-1}-3\alpha_sK_{\mathrm{eff}}(T^{k,m-1}-T_0)
    .
    \label{eq:volstrit}
\end{aligned}
\end{equation}
Then we have
\begin{equation}
    K_{\mathrm{eff}}\eps_\mathrm{vol}(\mbf{u}^{k,m}) = \frac{\alpha}{K_{\mathrm{eff}}}( p^{k,m} - p^{k,m-1}) + 3\alpha_s(T^{k,m} - T^{k,m-1}).
    \label{eq:volstrit_1}
\end{equation}
Substituting Eq.~\eqref{eq:volstrit_1} into Eq.~\eqref{eq:weak_form_of_fluid_flow_model_stagg}, the fluid flow equation now depends only on $p^{k,m}$ as
\begin{equation}
    \begin{aligned}
        &\int_{\mit\Omega} \left(\frac{\alpha^2}{K}\frac{p^{k,m}-p^{k,m-1}}{\Delta t}+3\alpha\alpha_s\frac{T^{k,m}-T^{k,m-1}}{\Delta t}\right)\mbf{w}_pdV \\
    &+\int_{\mit\Omega}\left(\frac{1}{M_p}\frac{p^{k,m}-p^{k-1}}{\Delta t}-\frac{1}{M_T}\frac{T^{k,m}-T^{k-1}}{\Delta t}\right)\mbf{w}_pdV\\
    &+ \int_{\mit\Omega}\frac{\mbf{K}}{\mu}\nabla p^{k,m}\cdot\nabla\mbf{w}_{p} dV =  \int_{\mit\Omega} \left(Q_{f}-\alpha \frac{\eps_\mathrm{vol}(\mbf{u}^{k,m-1})-\eps_\mathrm{vol}(\mbf{u}^{k-1})}{\Delta t}\right)\mbf{w}_{p} dV-\int_{\partial_{N}\mit\Omega}q_{n}\mbf{w}_{p} dS.
        \label{eq:eq:weak_form_fix_of_fluid_flow_model_stagg}
    \end{aligned}
\end{equation}
\vspace{20pt}

$\bullet$\:$\bm{\textbf{Step}\:4}$: Given\: $(T^{k,m}, p^{k,m}, \upsilon^{k,m})$, solve for\: $\mbf{u}^{k,m}$:
\begin{equation}
    \label{weak_form_dis_stagg}
    \int_{\mit\Omega}\nabla\mbf{w}_u\cdot\left[g(\upsilon^{k,m})\mbf{\sigma}_e^{k,m} -\alpha(\upsilon^{k,m})p^{k,m}\mbf{I}\right]\, \mathd V - \int_{\mit\Omega}\mbf{b}\cdot\mbf{w}_u\, \mathd V - \int_{\mathcal{C}_N}\bar{\mbf{t}}\cdot\mbf{w}_u\, \mathd S = \mbf{0}.
\end{equation}

The proposed thermo-hydro-mechanical phase-field model has been implemented in an open source code, OpenGeoSys~\cite{bilke_2023_7716938}, and discretization details are provided in \ref{B}.

\subsection{Isotropic diffusion stabilization method}
For hydraulic fracturing problems, a mass source term imposes a strong advective term in heat transfer, and the standard Galerkin finite element method may generate non-symmetric coefficient matrices, which lead to numerical instabilities and spurious oscillations~\cite{tezduyar2002calculation,donea2003finite}.
To address this issue, Burman and Ern (2002) proposed the isotropic diffusion method~\cite{burman2002nonlinear}, which adds an artificial isotropic balancing dissipation to the diffusion coefficient and forces the P\'eclet number to be smaller than 1. 
The isotropic balancing dissipation is defined as
\begin{equation}
    \mathbf{K}_{s} = \frac{1}{2} s ||\mbf{q}||h_e \mathbf I
    ,
\end{equation}
where $s \in [0,1] $ is the tuning parameter and is set to $s = 0.15$ in this study.

%% file: verif.tex
This section aims to verify our numerical model against known closed-form solutions. 
No closed-form solution is known to couple all thermo-hydro-mechanical components with fracture propagation.
Thus, we present three verification tests that test different aspects of the model: (1) Terzaghi's problem for the hydro-mechanical part, (2) the thermal consolidation problem for thermo-hydro-mechanical part, and (3) the KGD problem for the hydraulic fracture propagation part.

\subsection{Terzaghi's consolidation problem}
\begin{figure}[H]
		\centerline{\includegraphics[scale=0.6]{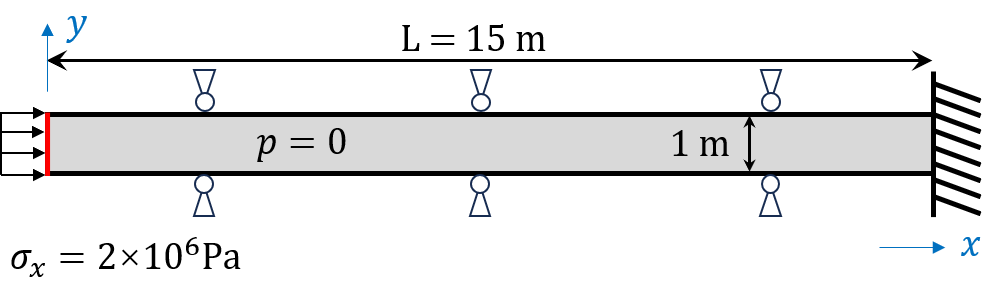}}
		\caption{\centering{Schematic of Terzaghi's consolidation problem.}}
		\label{fig:terzaghi_sche}
\end{figure}

Firstly, we verified the hydro-mechanical module (i.e.~the thermal and fracture modules were turned off) against Terzaghi's consolidation problem. 
Fig.~\ref{fig:terzaghi_sche} illustrates a fluid-saturated soil column's geometry and boundary conditions. 
For the mechanical boundary conditions, a constant stress $\sigma_x = 2$ MPa was applied on the left boundary while the normal displacement was fixed on the right boundary. 
For the flow boundary conditions, all the boundaries were set to be no-flow except the left edge, where the pressure was set to 0 MPa (drained condition). 
Table~\ref{Terzaghi_para} lists the material parameters. 
The analytical solutions for pressure and displacement profile evolution are given as~\cite{biot1941general}
\begin{flalign}
    & \ p(x,t)=\frac{4d\sigma_x}\pi\sum_{m=0}^\infty\left\{\frac{1}{2m+1}exp\left(-\frac{(2m+1)^2\pi^2}{4L^2}ct\right)sin\left(\frac{(2m+1)\pi x}{2L}\right)\right\}&\nonumber
\end{flalign}

\begin{flalign}
    & \ u(x,t)=c_md\sigma_k\left[L-x-\frac{8L}{\pi^2}\sum_{m=0}^\infty\left\{\frac{1}{(2m+1)^2}exp\left(-\frac{(2m+1)^2\pi^2}{4L^2}ct\right)\cos\left(\frac{(2m+1)\pi x}{2L}\right)\right\}\right]+b\sigma_x(L-x)&\nonumber
\end{flalign}
where
\begin{equation}
    a=\frac{(1+v)(1-2v)}{E(1-v)}, \quad S=\frac{\alpha-\phi}{K_{s}}+\phi_m c_{f}, \quad b=\frac{a}{1+a\alpha^{2}/S}, \quad d=\frac{a-b}{a\alpha}, \quad c=\frac{K_{m}}{(a\alpha^{2}+S)\mu}, \quad c_{m}=\frac{a-b}{d}.\nonumber
\end{equation}

\begin{figure}[H]
		\centerline{\includegraphics[scale=0.5]{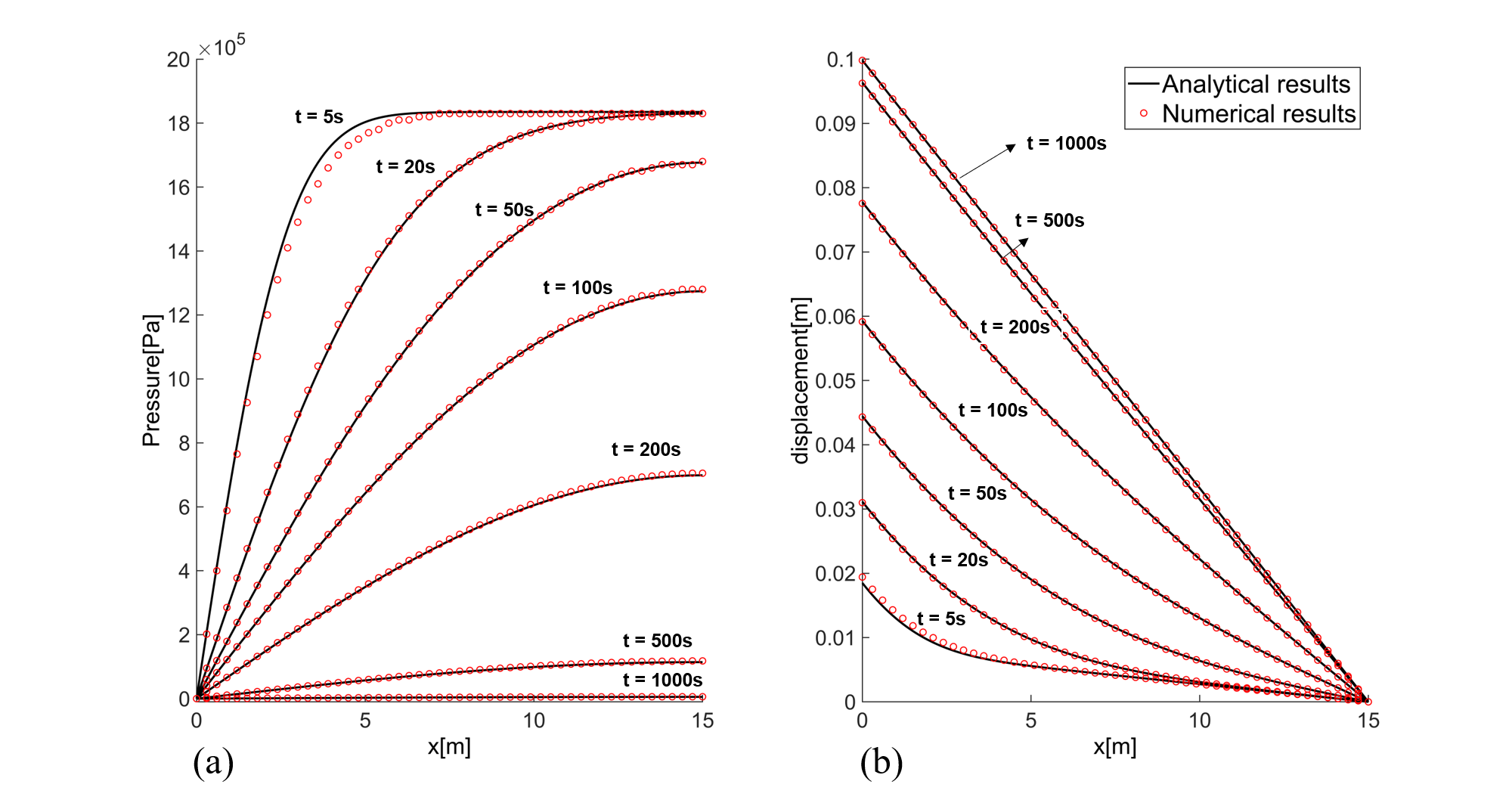}}
		\caption{Comparisons of numerical results and analytical solution of Terzaghi’s consolidation problem for (a) pressure and (b) displacement.}
		\label{fig:terzaghi_result}
\end{figure}

The total simulation time was 1000~s with a time increment of 1~s.
Simulated pressure and displacement closely match the analytical solution (Fig.~\ref{fig:terzaghi_result}).

\begin{table}[!ht]
    \renewcommand{\arraystretch}{1.5}
    \setlength{\belowcaptionskip}{0.5cm}
	\centering
	\caption{\centering{Parameters for Terzaghi's problem.}}
    \label{Terzaghi_para}
	\begin{tabular}{p{6cm}p{4cm}l}
		\hline
		Input parameters & Value & Unit \\ 
		\hline
		Young's modulus ($E$) & $0.3$	& GPa \\		 
		Poisson's ratio ($v$)& 0 & - \\
		Biot coefficient ($\alpha_m$) & 1	& - \\
		Porosity ($\phi_m$) & 0.3	& - \\
        Permeability ($\mathbf{K}_m$) & $2\times 10^{-12}$	& m$^2$ \\
		\hline
	\end{tabular}
\end{table}
\subsection{Thermal consolidation problem}
\begin{figure}[H]
		\centerline{\includegraphics[scale=0.7]{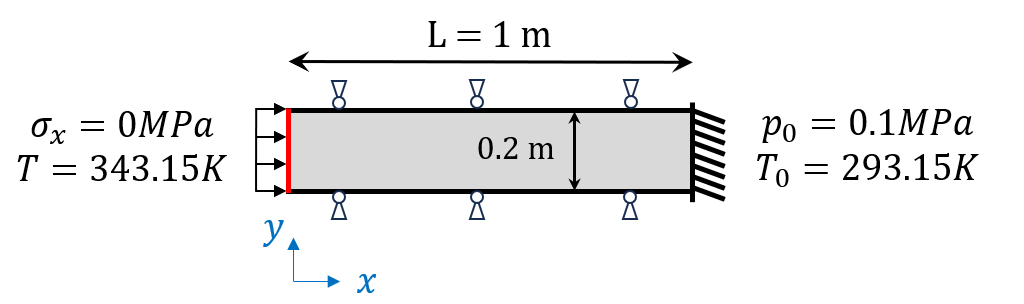}}
		\caption{\centering{Schematic of the thermal consolidation problem.}}
		\label{fig:thermal_consolidation_sche}
\end{figure}

In this benchmark example, we verified the thermo-hydro-mechanical coupling module against the analytical solution provided by~\cite{shi2019numerical} for the thermal consolidation problem~\cite{bing2005one}. 
Fig.~\ref{fig:thermal_consolidation_sche} shows a consolidated soil column $[ 0 \, \mrm{m},\, 1 \, \mrm{m} ] \times [ 0 \, \mrm{m},\, 0.2 \, \mrm{m} ]$ with the initial temperature of 293.15 K and the initial pressure of 0.1~MPa. 
A constant temperature of 343.15~K was applied on the left edge while the pressure was fixed as 0~MPa.
The remaining three boundaries were considered impervious and insulated.
For the mechanical boundary conditions, the displacements were constrained in the y-direction at $y=0$ m and $y=0.2$ m while the displacements were constrained in both the x- and y-directions at $x=1$ m.
Table~\ref{thermal_consolidation_para} lists the properties of soil and fluid. 

\begin{figure}[H]
		\centerline{\includegraphics[scale=0.55]{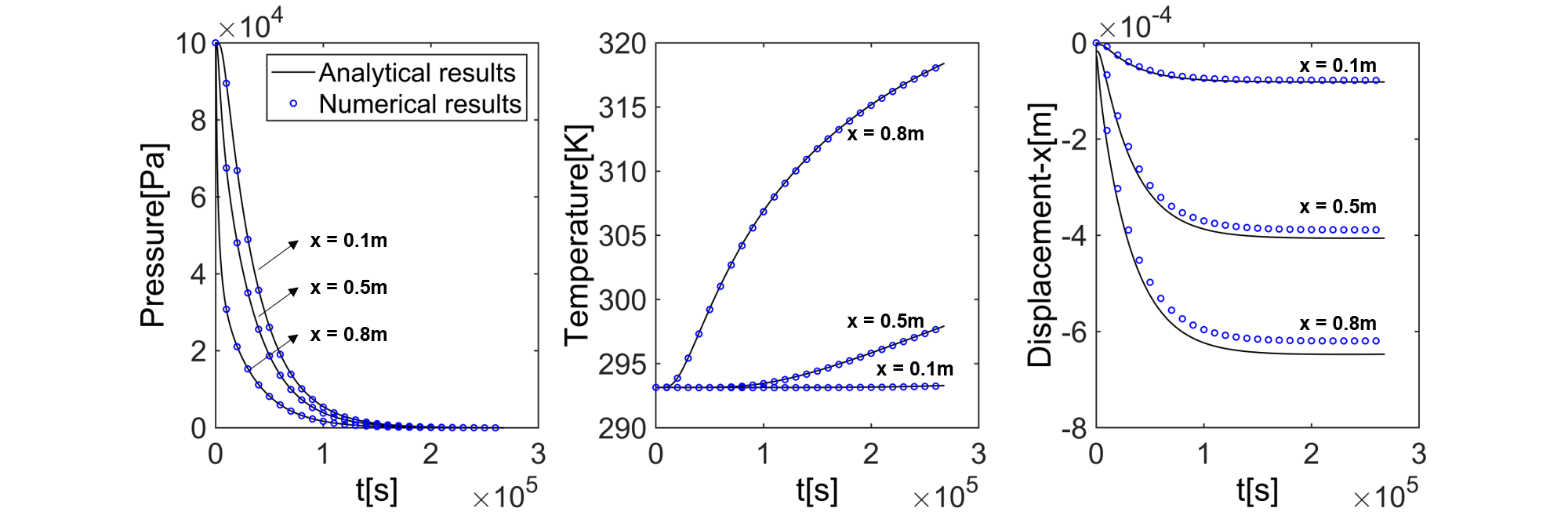}}
		\caption{\centering{Comparison of numerical results and analytical solution of Terzaghi’s problem for (a) pressure, (b) temperature, and (c) displacement.}}
		\label{fig:thermal_consolidation_res}
\end{figure}

Fig.~\ref{fig:thermal_consolidation_res} compares the simulated evolution of pressure, temperature, and displacement at different locations against the analytical solution.
Although the displacement evolution shows slight differences at a distance, the results agree well.
\begin{table}[!ht]
    \setlength{\belowcaptionskip}{0.5cm}
    \renewcommand{\arraystretch}{1.5}
	\centering
	\caption{\centering{Parameters for the thermal consolidation problem.}}
    \label{thermal_consolidation_para}
	\begin{tabular}{p{9cm}p{4cm}l}
		\hline
		Input parameters & Value & Unit \\ 
		\hline
		Young's modulus ($E$) & 60 & MPa \\		 
		Poisson's ratio ($v$)& 0.4 & - \\
		Biot coefficient ($\alpha_m$) & 1	& - \\
		Porosity ($\phi_m$) & 0.4	& - \\
        Permeability ($\mathbf{K}_m$) & 1e-16	& m$^2$ \\
        Thermal conductivity ($\lambda$) & 0.5 & W/(m$\cdot$ K) \\
        Thermal expansivity of soil ($\alpha_s$) & 3e-7	& $-$ \\
        Specific heat capacity of soil and fluid ($c_{p,s}, c_{p,f}$) & 800, 4200 & J/(kg$\cdot$ K) \\
		\hline
	\end{tabular}
\end{table}

\subsection{Fluid-driven fracture propagation in plane-strain}
\begin{figure}[H]
		\centerline{\includegraphics[scale=0.7]{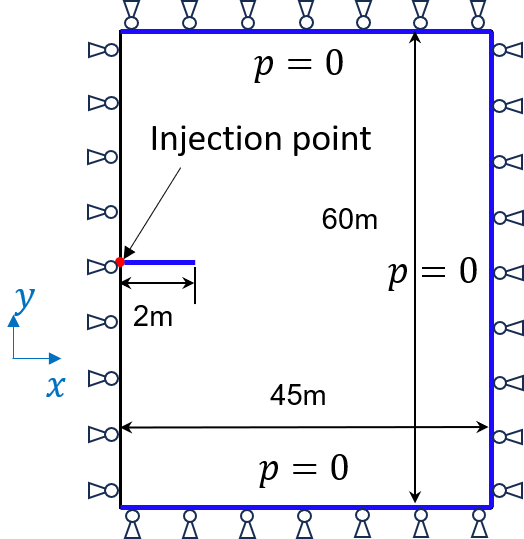}}
		\caption{Geometry and boundary conditions of the KGD hydraulic fracture problem with a half symmetry across the $y$-axis.}
		\label{fig:KGD}
\end{figure}

As our last benchmark, we verified our model's hydro-mechanical module with a propagating fracture (phase-field) in the plane strain condition against a so-called KGD (Kristianovich–Geertsma–de Klerk) model without fluid leak-off.
Assuming a symmetry over the $y$-axis, a line fracture $[ 0 \, \mrm{m},\, 2 \, \mrm{m} ] \times \{ 30 \, \mrm{m} \}$ is considered in a domain $[ 0 \, \mrm{m},\, 45 \, \mrm{m} ] \times [ 0 \, \mrm{m},\, 60 \, \mrm{m} ]$ (Fig.~\ref{fig:KGD}).
We considered incompressible fluid ($c_f = 0$) injection into the impermeable elastic medium ($\alpha_m = 0$ and $\phi_m = 0$) to propagate the line fracture.
The smallest element size is 0.05 m~and $\ell / h_e = 4$. 
In the first 10 time steps, $\Delta t = 0.01$~s and then $\Delta t = 0.1$~s in the remaining simulation time. 
Table~\ref{KGD_para} lists the mechanical and flow parameters. 

\begin{table}[!ht]
    \renewcommand{\arraystretch}{1.5}
    \setlength{\abovecaptionskip}{0cm}
    \setlength{\belowcaptionskip}{0.5cm}
	\centering
	\caption{\centering{Mechanical and flow parameters for KGD problem.}}
    \label{KGD_para}
	\begin{tabular}{p{9cm}p{4cm}l}
		\hline
		Input parameters & Value & Unit \\ 
		\hline
		Young's modulus ($E$) & 17 & GPa \\		 
		Poisson's ratio ($v$)& 0.2 & - \\
        Permeability ($\mathbf{K}_m$) & $1 \times 10^{-18}$	& m$^2$ \\
        Fluid viscosity ($\mu$) & $1 \times 10^{-8}$ & Pa$\cdot$s \\
        Injection rate ($Q$) & $2 \times 10^{-3}$	& m$^2$/s\\
        Critical surface energy release rate (G$_c$) & 300 & N/m \\
		\hline
	\end{tabular}
\end{table}

The hydraulic fracture propagation in this setting is considered a toughness-dominated regime~\cite{detournay2003near} in which the fluid viscous dissipation is negligible compared to the energy release by fracture propagation.
To judge the fracture propagation regime, we can use the dimensionless viscosity $\mathcal{M}$ defined for the KGD fracture as~\cite{garagash2006plane}
\begin{equation}
    \mathcal{M}=\frac{\mu^{\prime}Q}{E^{\prime}}\left(\frac{E^{\prime}}{K^{\prime}}\right)^{4}
    \label{dimensionless_M}
    ,
\end{equation}
with $K^{\prime} = \sqrt{\frac{32G_{\mathrm{c}} E^{\prime}}\pi}$, $\mu^{\prime} = 12\mu$, $E^{\prime} = \frac{E}{1 - \nu^2}$.
And if $\mathcal{M}  < \mathcal{M}_c = 3.4\times 10^{-3}$,  the KGD fracture is toughness dominated.
In our setting, $\mathcal{M} = 3.8\times 10^{-7}$ and $\mathcal{M}  < \mathcal{M}_c$\footnote{The effective critical surface energy release rate $G_c$ in the phase-field model needs to be adjusted to account for the discretization~\cite{bourdin2008variational, tanne2018crack, yoshioka2020crack}. With the mesh size $h_e$ and the characteristic length $\ell$, the modification is given by
\begin{equation}
   G_{\mathrm{c}}^{\mathrm{eff}}=G_{\mathrm{c}}\left(1+\frac{h_e}{4c_{n}\ell}\right). \nonumber
\end{equation}
}.



\begin{figure}[H]
    \centering
      \subfloat[Pressure response at the injection point.]
      {\includegraphics[scale=0.5]{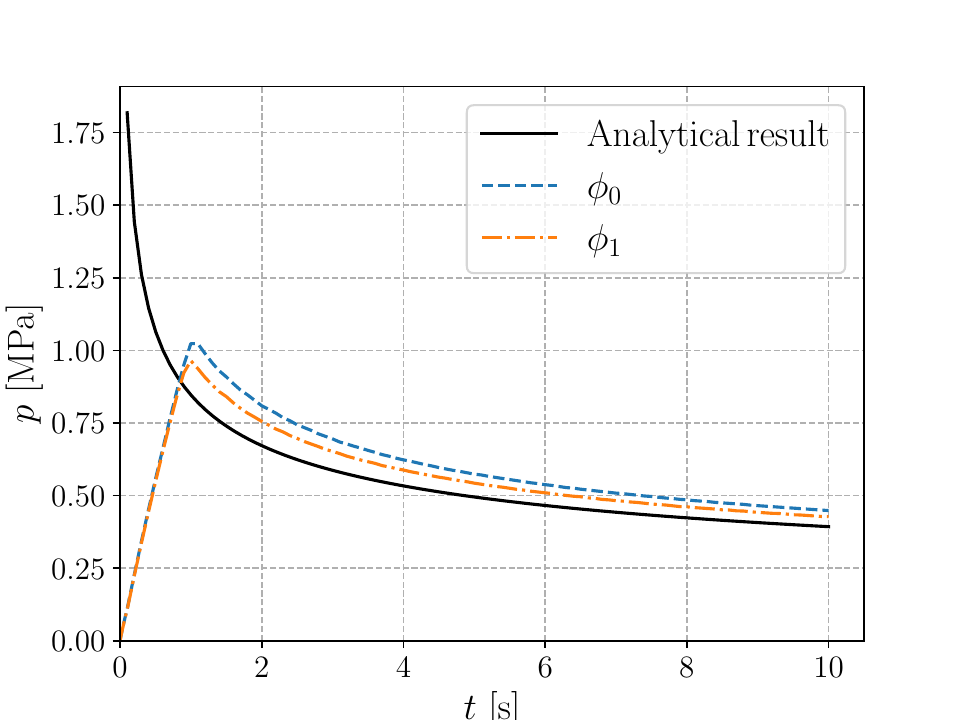}\label{fig:KGD_p}}
      \subfloat[Fracture length evolution.]
      {\includegraphics[scale=0.5]{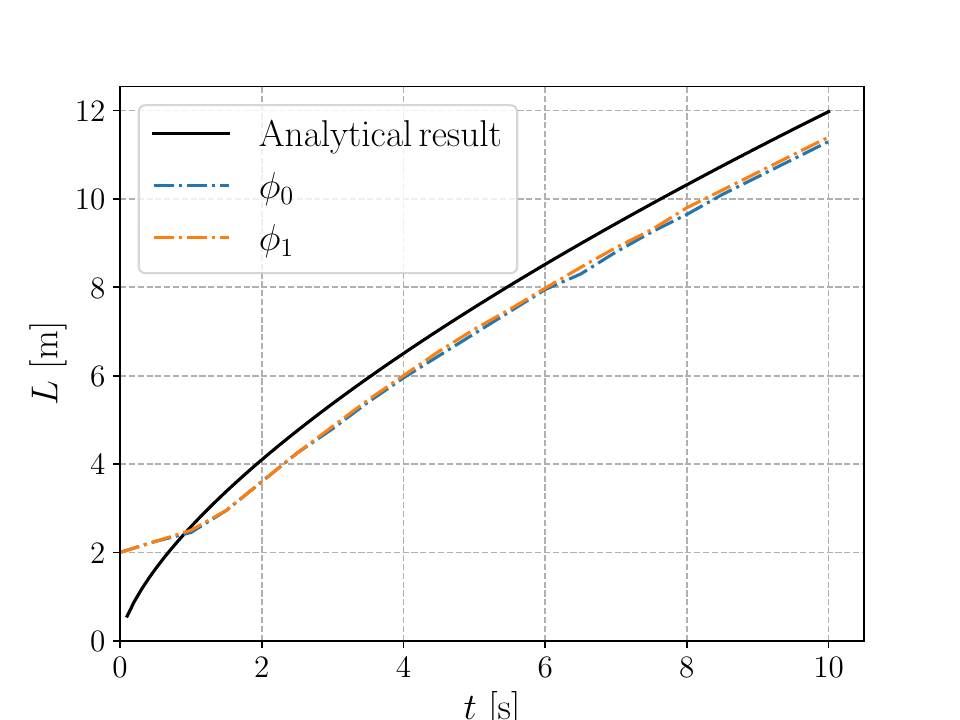}\label{fig:KGD_l}}
      \quad
      \subfloat[Fracture width evolution at the injection point.]
      {\includegraphics[scale=0.5]{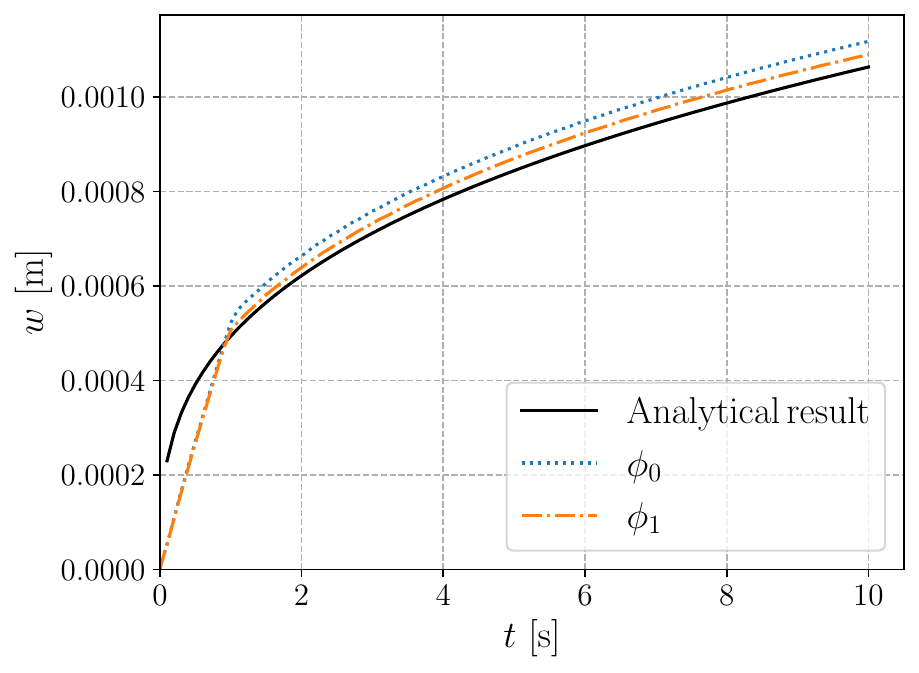}\label{fig:KGD_w}}
      \caption{Comparisons of the KGD hydraulic fracture model between different degraded porosity models.}
      \label{fig:Compare_results_KGD}
\end{figure}

Here, we compare the simulation results from 2 different porosity models -- $\phi_0$ (Eq.~\ref{eq:pf_porosity_0}) and $\phi_1$ (Eq.~\ref{eq:pf_porosity}) -- against the analytical solution~\cite{garagash2006plane} in Fig.~\ref{fig:Compare_results_KGD}.
The results show that our proposed porosity model ($\phi_1$) improves the solution accuracy in the pressure and fracture width at the injection point and the fracture length evolution.

As theoretically demonstrated in~\cite{bourdin2000numerical}, the phase-field approximation approaches the sharp crack representation as $\ell \rightarrow 0$.
In practice, however, $h_e$ has to be greater than $\ell$ (i.e. $\ell/h_e > 1$)~\cite{bourdin2012variational}.
This implies that simulation results converge to the sharp crack result (analytical solution) as we refine the mesh ($h_e$), and it was also demonstrated in~\cite{YOU2023116305}. 
For the same element resolution, Fig.~\ref{fig:Compare_results_KGD} suggests that our proposed porosity model ($\phi_1$) represents the sharp fracture behaviors more closely than the phase-field dependent porosity model ($\phi_0$). 


%% file: numexperi.tex
This section presents several numerical fracturing examples to study the complex interactions between thermal, hydraulic, and mechanical effects, including an advection-dominated problem and propagation of a single fracture under different injection temperatures with and without a weak interface. 

\subsection{Stabilization for the advection-dominated problem}
\begin{figure}[H]
    \centerline{\includegraphics[scale=0.6]{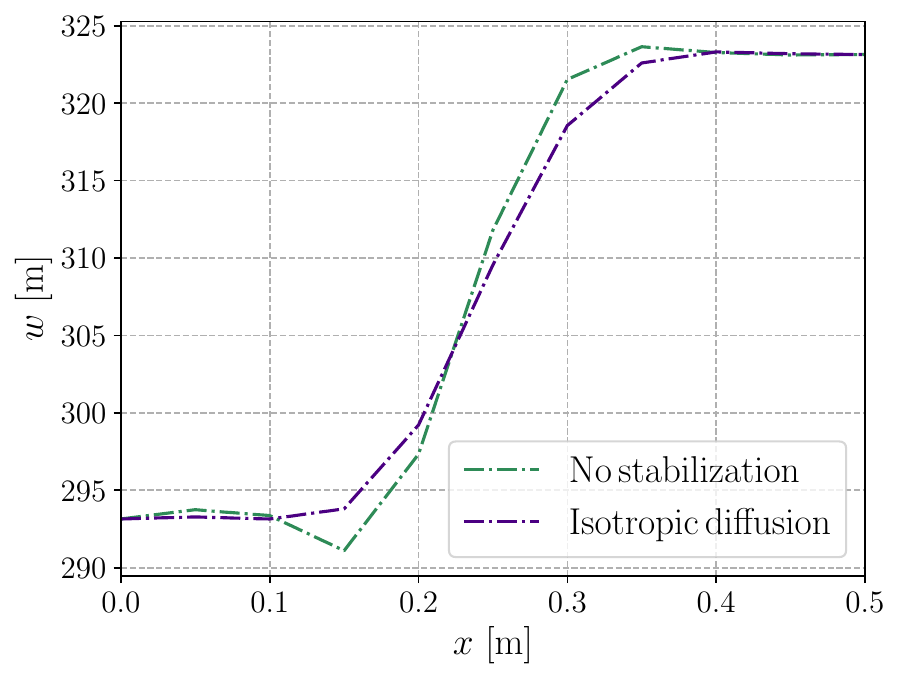}}
    \caption{Temperature profiles along the fracture with and without the isotropic diffusion stabilization method.}
    \label{fig:T_stab}
\end{figure}

First, we tested the isotropic diffusion stabilization method for an advection-dominated problem. 
With the same parameters in the KGD verification example (Table~\ref{KGD_para}), we simulated cold fluid injection with a temperature difference of 30~K and compared the temperature profiles along the propagating fracture with and without the stabilization method (Fig.~\ref{fig:T_stab}). 
Without the stabilization term, the temperature profile shows unrealistic behaviors such as dropping below the injection temperature or increasing above the initial reservoir temperature.
On the other hand, with the stabilization term, the model was able to remove such spurious temperature behaviors, generating a smooth temperature profile along the fracture.

\subsection{Cold fluid injection into a single fracture}
\label{sec:single-cracked model}
\begin{figure}[H]
    \centerline{\includegraphics[scale=0.7]{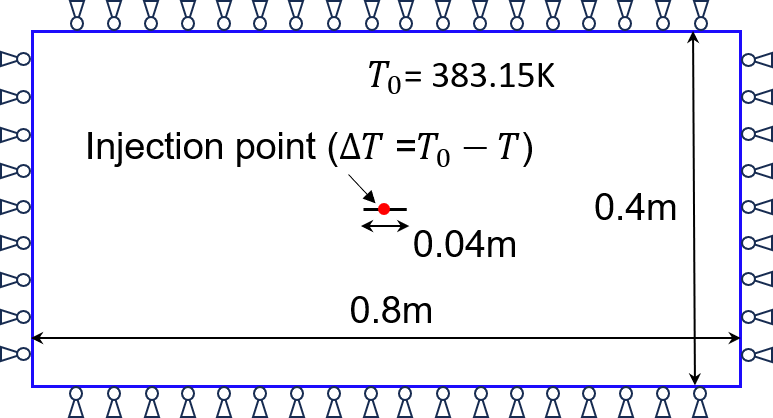}}
    \caption{\centering{The schematic of a single fracture in the middle of the domain.}}
    \label{fig:single_cracked_model}
\end{figure}

In this example, we simulated cold fluid injection into a single fracture in a poro-elastic medium as illustrated in Fig.~\ref{fig:single_cracked_model}. 
Consider a computation domain $[ 0 \, \mrm{m},\, 0.8 \, \mrm{m} ] \times [ 0 \, \mrm{m},\, 0.4 \, \mrm{m} ]$ with an initial fracture $[ 0.38 \, \mrm{m},\, 0.42 \, \mrm{m} ] \times \{ 0.4 \, \mrm{m} \}$.
With fluid injection at a rate of $2 \times 10^{-5}$~m$^2$/s at $( 0.4 \, \mrm{m},\, 0.2 \, \mrm{m} )$, we imposed 4 different temperature differences ($\Delta T=$ 0~K, 30~K, 60~K and 90~K) between the injection fluid ($T$) and the initial reservoir ($T_0 = 383.15$~K) by decreasing the injection temperature. 
The initial pressure in the domain is 0~MPa, and all the boundaries are in the drained condition ($p=0$~MPa) with the normal displacement constrained. 
The material properties are listed in Table~\ref{single_cracked_para}. 

\begin{table}[!ht]
    \renewcommand{\arraystretch}{1.5}
    \setlength{\abovecaptionskip}{0cm}
    \setlength{\belowcaptionskip}{0.5cm}
	\centering
	\caption{\centering{Parameters for single-cracked model.}}
    \label{single_cracked_para}
	\begin{tabular}{p{9cm}p{4cm}l}
		\hline
		Input parameters & Value & Unit \\ 
		\hline
		Young's modulus ($E$) & 17 & GPa \\		 
		Poisson's ratio ($v$)& 0.2 & - \\
        Biot's coefficient ($\alpha_m$)& 0.6 & - \\
        Permeability ($\mathbf{K}_m$) & 1e-16	& m$^2$ \\
        Thermal conductivity of soil and fluid ($\lambda$) & 3, 0.5	& W/(m$\cdot$K) \\
        Thermal expansivity ($\alpha_s$) & $8\times 10^{-6}$	& - \\
        Specific heat capacity of soil and fluid ($c_{p,s}$, $c_{p,f}$) & 800, 4200& J/(kg$\cdot$K) \\
        Fluid viscosity ($\mu$) & $1 \times 10^{-4}$& Pa$\cdot$s \\
        Critical surface energy release rate ($G_c$) & 100 & N/m \\
		\hline
	\end{tabular}
\end{table}

\begin{figure}[H]
    \centering
      \subfloat[Pressure response at the injection point.]
      {\includegraphics[scale=0.5]{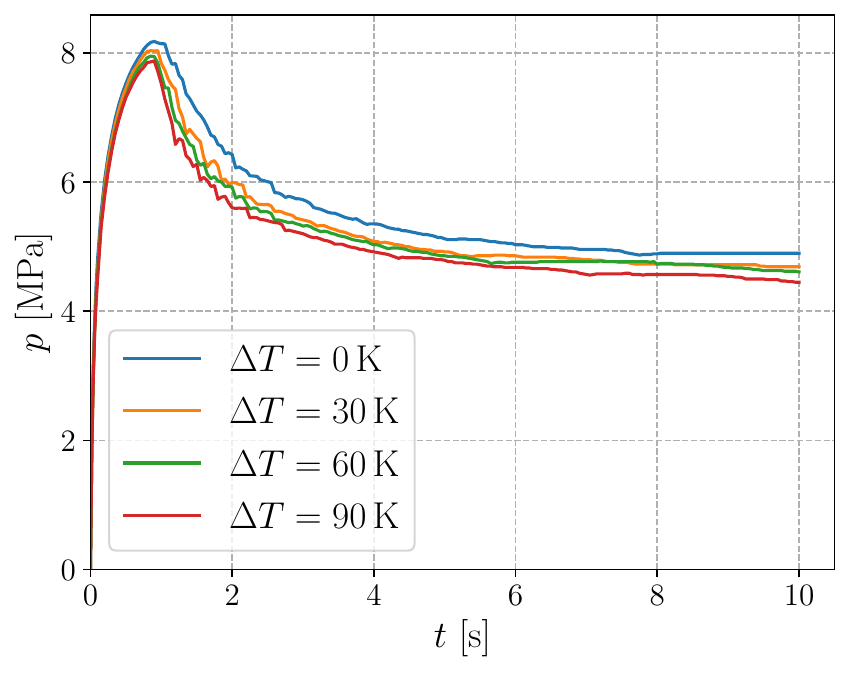}\label{fig:single_frac_p}}
      \subfloat[Fracture length evolution.]
      {\includegraphics[scale=0.5]{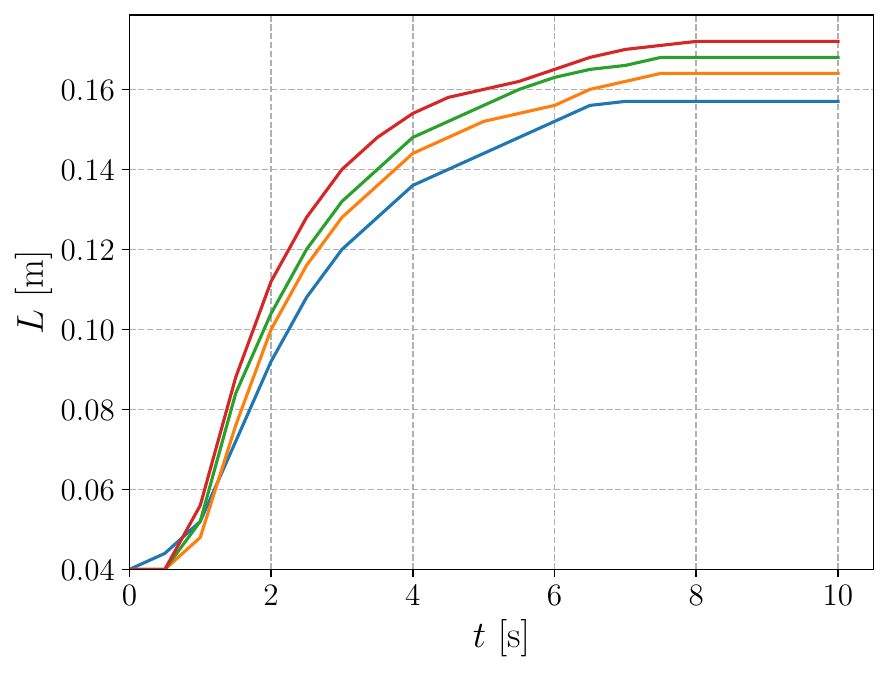}\label{fig:single_frac_l}}
      \quad
      \subfloat[Fracture width evolution at the injection point.]
      {\includegraphics[scale=0.5]{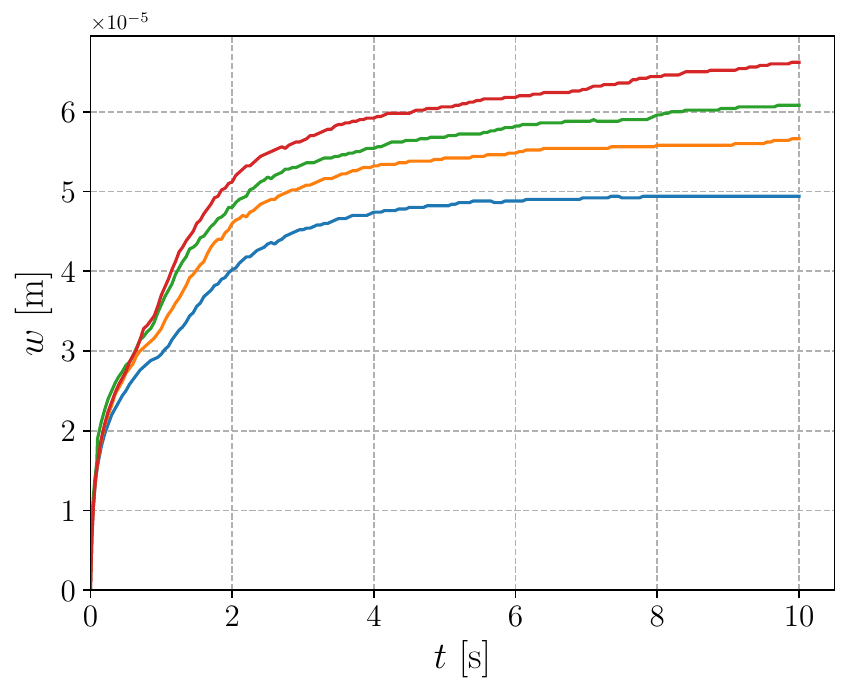}\label{fig:single_frac_w}}
      \caption{Propagation of a single fracture under different temperature differences ($\Delta T =$ 0~K, 30~K, 60~K and 90~K).}
      \label{fig:Compare_results_single_frac}
\end{figure}

As the injection temperature decreases, the critical pressure for fracture propagation decreases slightly (Fig.~\ref{fig:Compare_results_single_frac}(a)). 
This is because the contraction of rock contributes to the strain energy, and less pore fluid energy is required to reach the critical pressure. 
Once the fracture starts to propagate, the pressure diffuses more quickly with lower injection temperature because of the longer fracture length and the larger crack opening as shown in Fig.~\ref{fig:Compare_results_single_frac}(b) and (c). 
The pressures oscillate after the onset of fracture propagation, compared to the no leak-off case in Fig.~\ref{fig:Compare_results_KGD}.
This is because the fracture pressure drops not only with fracture growth but also with fluid leak-off to the formation.  
As a result, the pressure drop is more unstable.

\begin{figure}[H]
    \centerline{\includegraphics[scale=0.5]{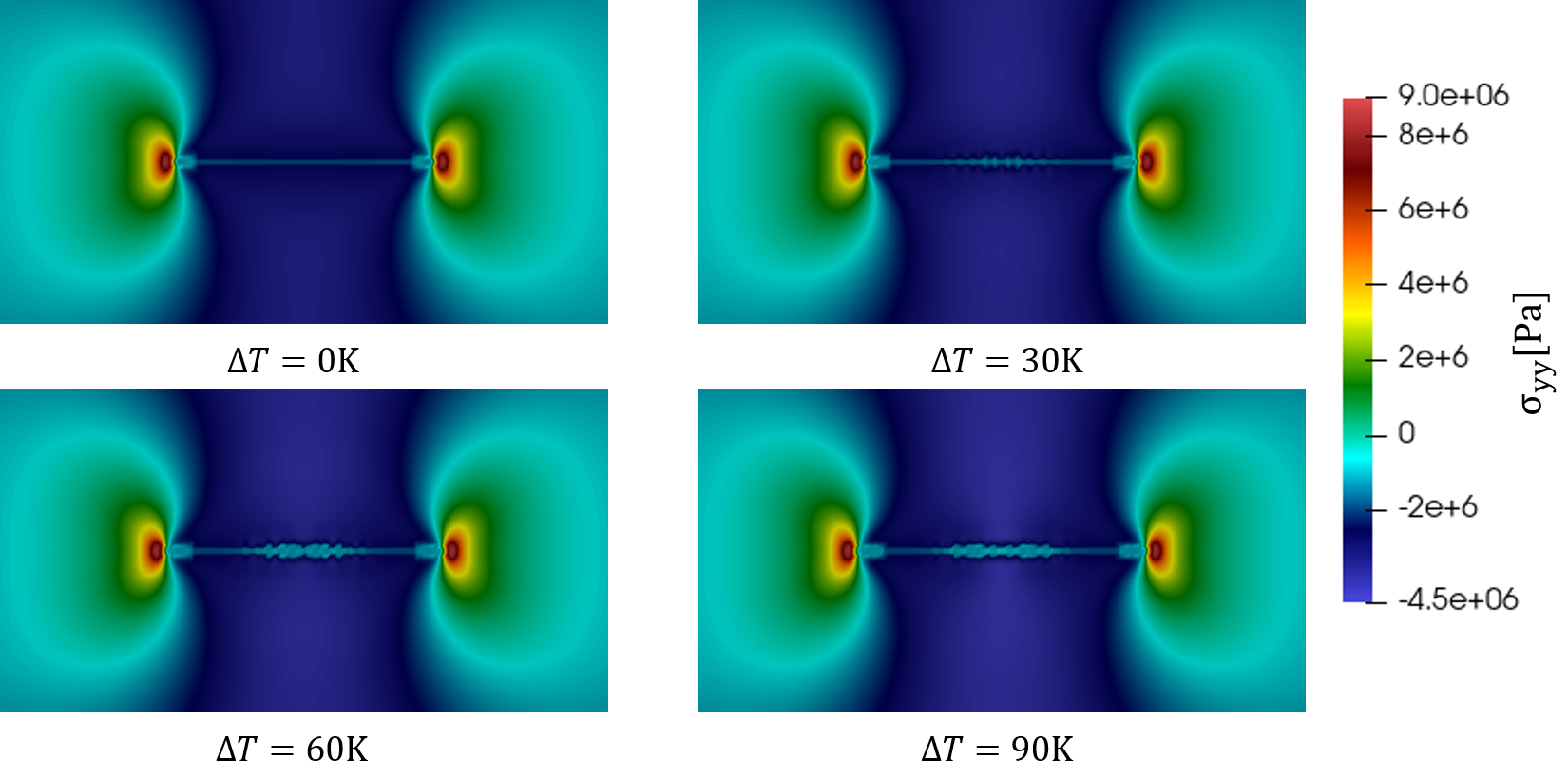}}
    \caption{The profiles of the stress in the $y$ component ($\sigma_{yy}$) with (a) $\Delta T = 0$ K (b) $\Delta T = 30$ K (c) $\Delta T = 60$ K, and (d) $\Delta T = 90$ K. }
    \label{fig:sigma_varia_tem}
\end{figure}

A lower injection temperature promotes the propagation of the fracture (Fig.~\ref{fig:Compare_results_single_frac}) as a temperature drop reduces the effective stress acting on the fracture boundaries. 
Fig.~\ref{fig:sigma_varia_tem} indicates how the effective stress in the $y$-direction is reduced around the injection point, which is more pronounced with a higher temperature difference.


\subsection{Interaction with an interface}
This section explores hydraulic fracture interactions with a weak interface (pre-existing natural fracture).
Consider the same computational domain in Sec.~\ref{sec:single-cracked model} with an initial fracture $[ 0.39 \, \mrm{m},\, 0.41 \, \mrm{m} ] \times \{ 0.4 \, \mrm{m} \}$.
We added a weak interface $\{ 0.04 \, \mrm{m} \} \times [ 0.15 \, \mrm{m},\, 0.25 \, \mrm{m} ] $ whose fracture toughness is lower than the intact rock ($G_c^\mrm{int}/G_c = 0.5$) following the interface modeling approach proposed in~\cite{yoshioka2021variational}.
The material properties and injection rate are the same as the previous case in Sec.~\ref{sec:single-cracked model}.

\begin{figure}[H]
    \centering
      \subfloat[$\Delta T = 0$~K.]
      {\includegraphics[scale=0.21]{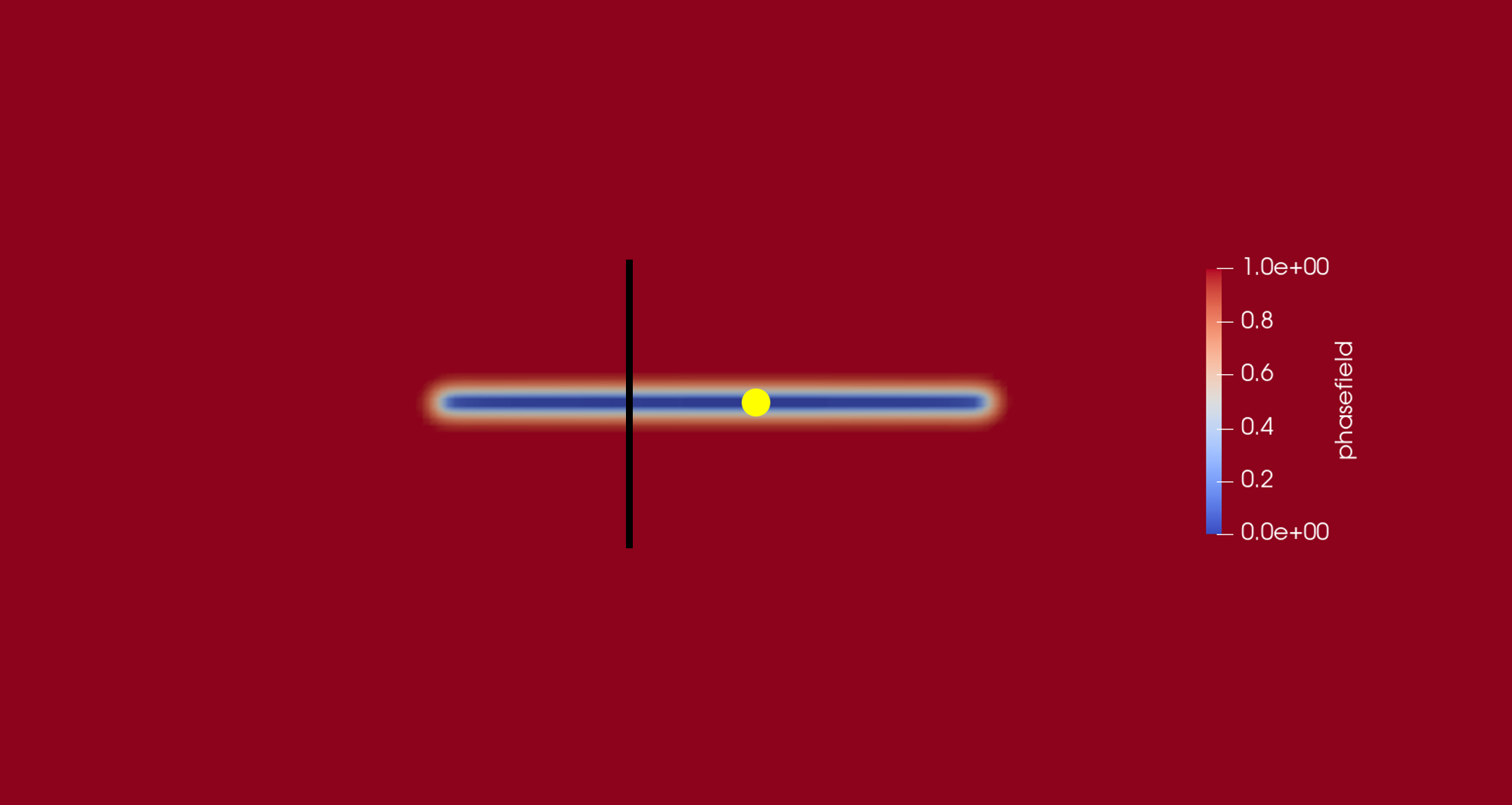}\label{fig:pf_nf_0K_inj}}
      ~
      \subfloat[$\Delta T = 30$~K.]
      {\includegraphics[scale=0.21]{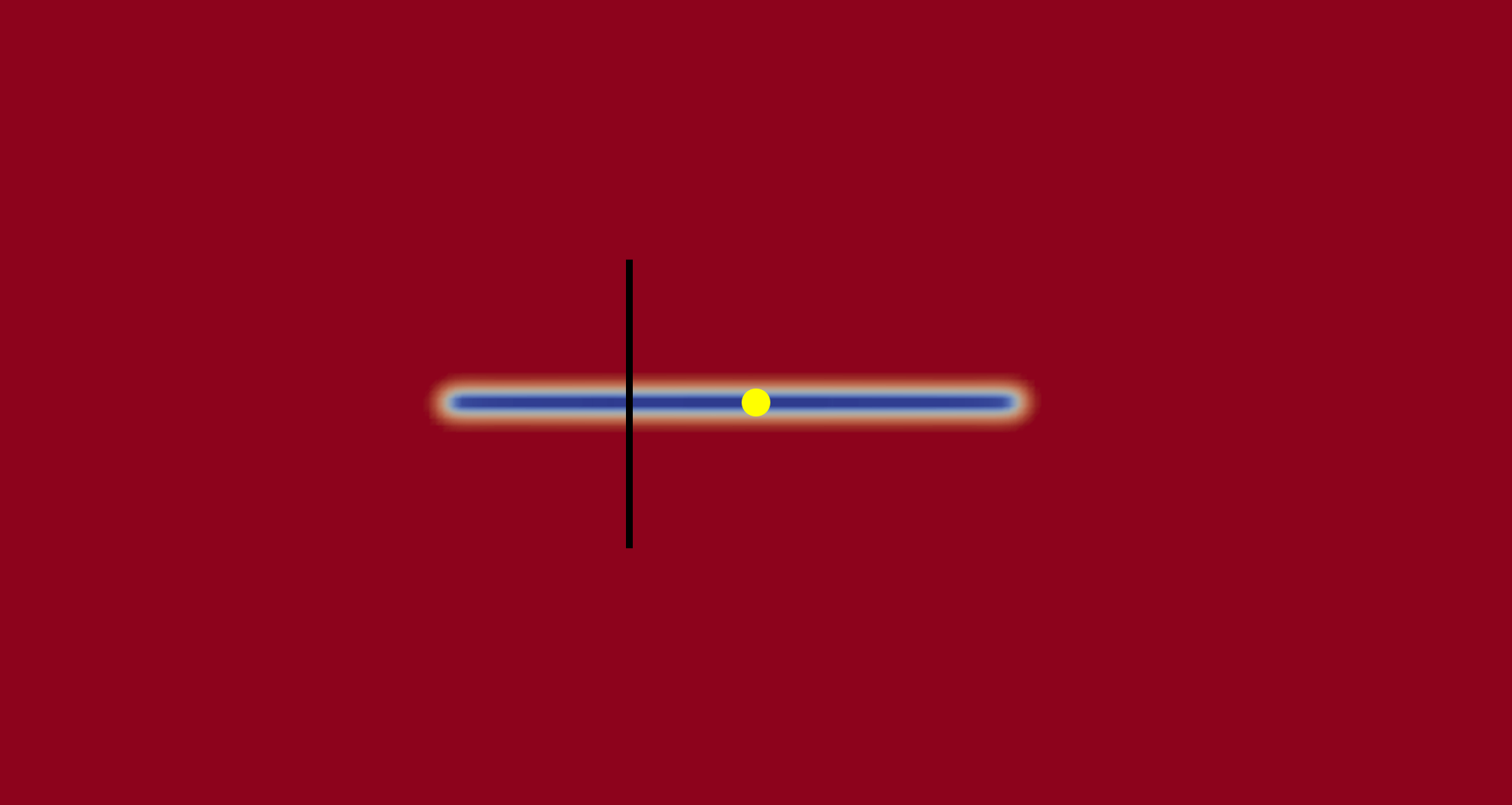}\label{fig:pf_nf_30K_inj}}
      \quad
      \subfloat[$\Delta T = 60$~K.]
      {\includegraphics[scale=0.21]{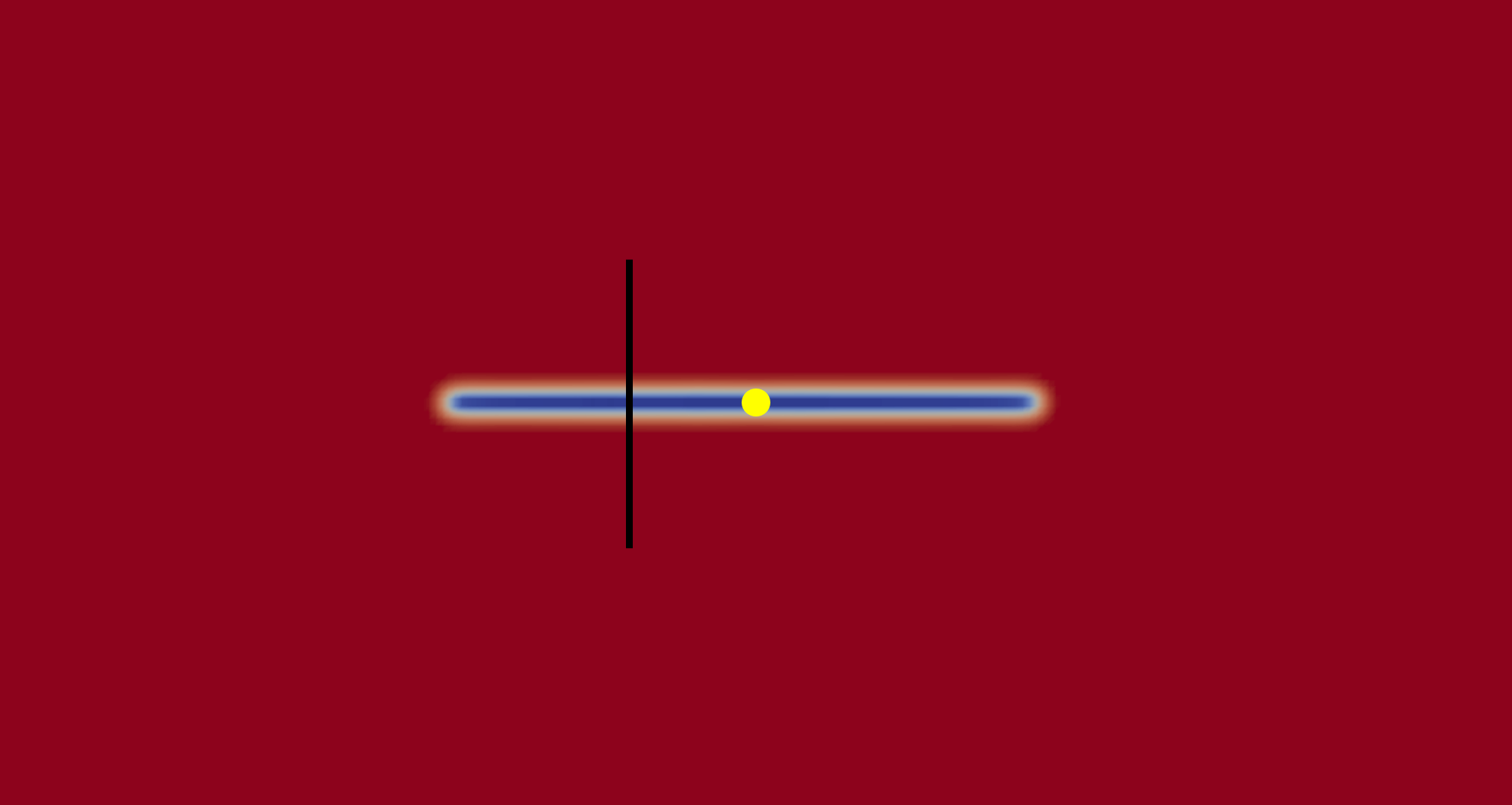}\label{fig:pf_nf_60K_inj}}
      ~
      \subfloat[$\Delta T = 90$~K.]
      {\includegraphics[scale=0.21]{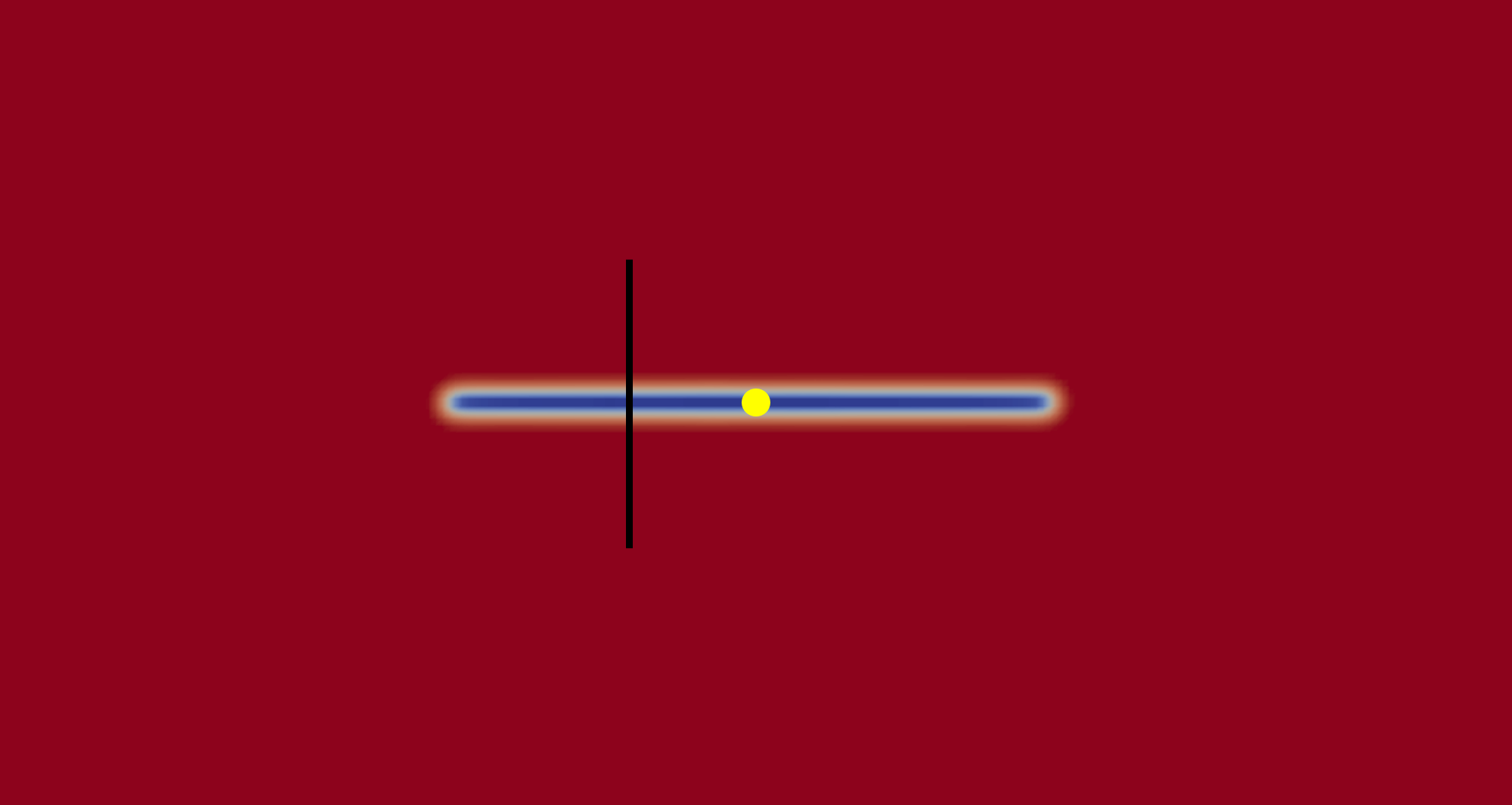}\label{fig:pf_nf_90K_inj}}
      \caption{{Hydraulic fracture interactions and propagation for (a) $\Delta T = 0$ K (b) $\Delta T = 30$ K (c) $\Delta T = 60$ K, and (d) $\Delta T = 90$ K. The injection point, located in the middle of the initial fracture, is marked with a yellow dot. In all the cases, the fracture propagation is attracted to the weak interface at time = 10 s.}}
      \label{fig:pf_nf_location}
\end{figure}

\begin{figure}[H]
    \centering
      \subfloat[Temperature difference is 0~K.]
      {\includegraphics[scale=0.5]{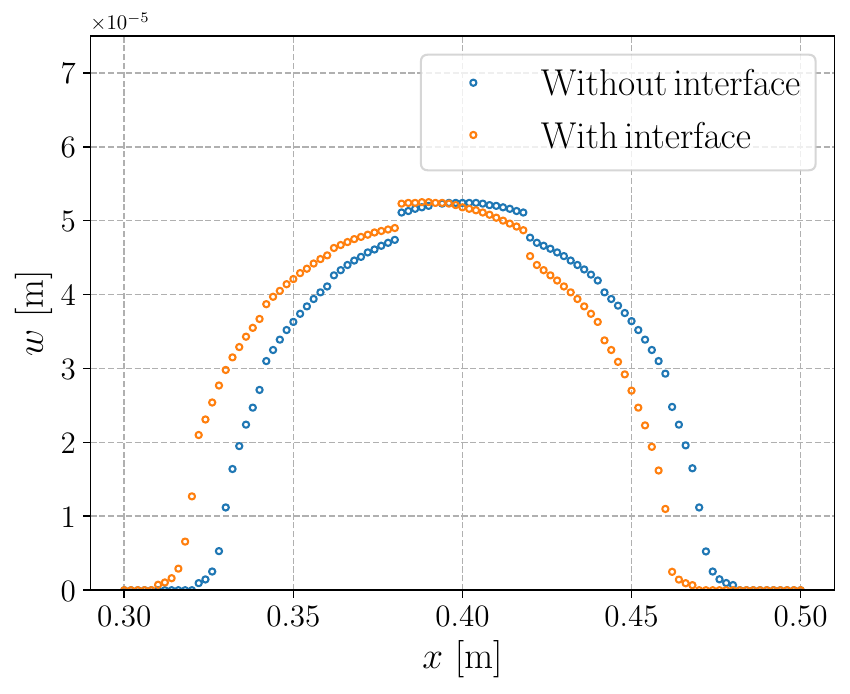}\label{fig:nf_0K}}
      \subfloat[Temperature difference is 30~K.]
      {\includegraphics[scale=0.5]{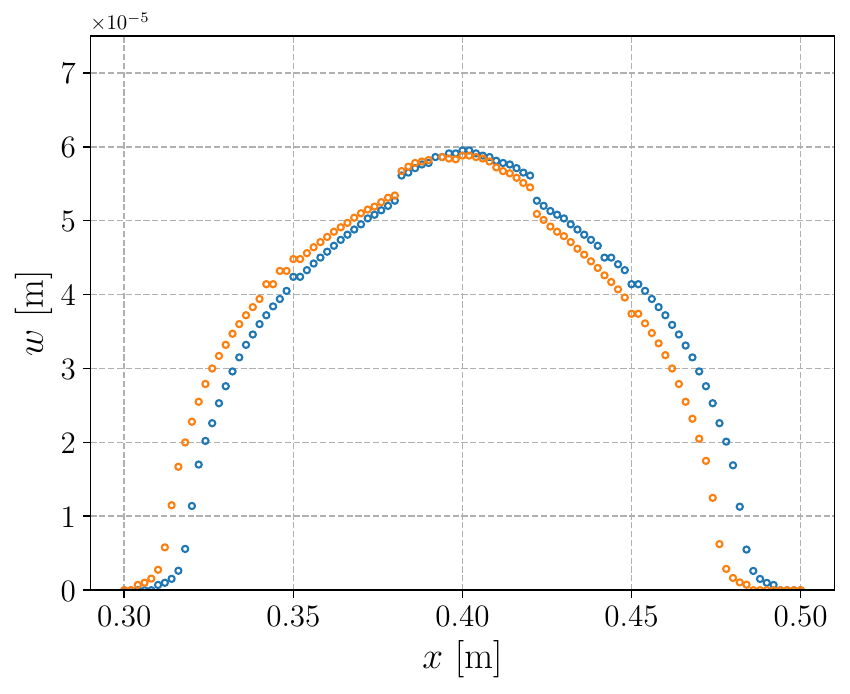}\label{fig:nf_30K}}
      \quad
      \subfloat[Temperature difference is 60~K.]
      {\includegraphics[scale=0.5]{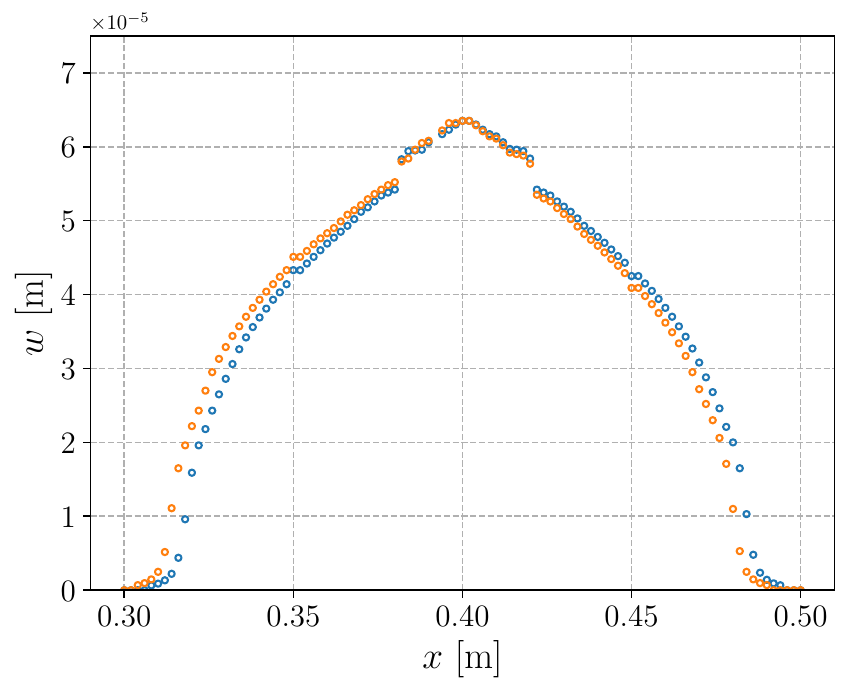}\label{fig:nf_60K}}
      \subfloat[Temperature difference is 90~K.]
      {\includegraphics[scale=0.5]{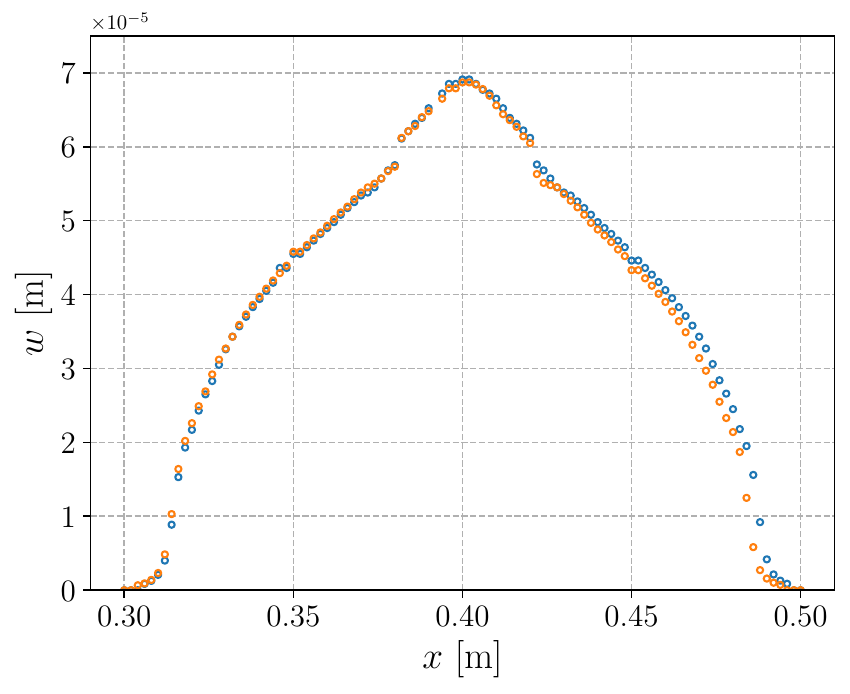}\label{fig:nf_90K}}
      \caption{{Comparisons of the fracture width ($\omega$) with and without the weak interface (natural fracture) for (a) $\Delta T = 0$ K (b) $\Delta T = 30$ K (c) $\Delta T = 60$ K, and (d) $\Delta T = 90$ K at time = 10 s.}}
      \label{fig:Compare_results_crack_opening}
\end{figure}

We simulated fracture propagation for 4 temperature differences (0~K, 30~K, 60~K, and 90~K) with and without the weak interface. 
With the weak interface, the fracture tends to grow more towards the interface in all the cases to some degree, though the weak interface does not alter the fracture propagation orientation (Fig.~\ref{fig:pf_nf_location}).
Fracture width comparisons at time = 10~s also show the attraction of the fractures to the weak interface (Fig.~\ref{fig:Compare_results_crack_opening}).

\begin{figure}[H]
    \centerline{\includegraphics[scale=0.7]{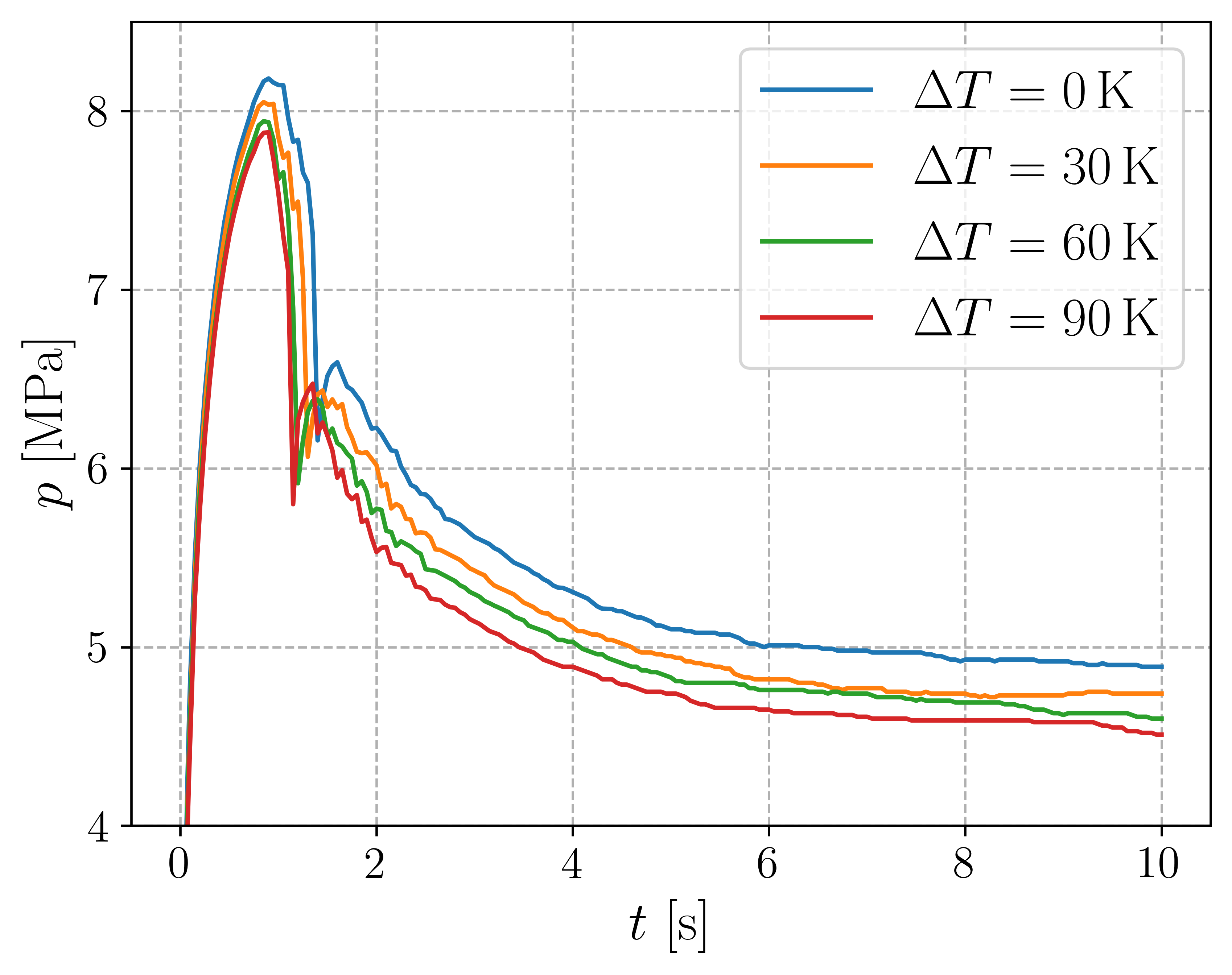}}
    \caption{Pressure responses at the injection point with the weak interface under different temperature differences (0~K, 30~K, 60~K, and 90~K).}
    \label{fig:p_inj}
\end{figure}

Though the fracture tends to grow towards the weak interface, colder injection temperatures (higher $\Delta T$s) inhibit this tendency.
As demonstrated in the single fracture case in the previous section (Fig.~\ref{fig:Compare_results_single_frac}), a higher $\Delta T$ leads to a longer and wider fracture with a lower propagation pressure.
Thus, with a higher $\Delta T$, the presence of the weak interface is overshadowed by the thermal stress impacts.
A higher temperature difference also induces larger fracture openings around the injection point where the temperature is colder, and the stress is reduced more (Fig.~\ref{fig:Compare_results_crack_opening}).
The reduced stress can also be seen in the pressure responses at the injection point (Fig.~\ref{fig:p_inj}), which show a pressure dip when the fracture reaches the interface. 
Similarly to the previous single fracture cases, the fracture propagation starts earlier with a higher temperature difference.


\begin{figure}[H]
    \centerline{\includegraphics[scale=0.7]{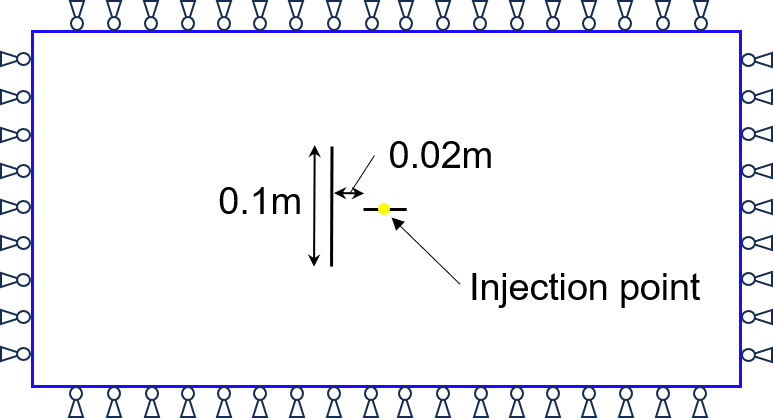}}
    \caption{The schematic of a single fracture in the middle with an orthogonal weak interface (natural fracture). The injection point is marked with a yellow dot.}
    \label{fig:single_cracked_model_nature_frac}
\end{figure}

%% file: conc.tex
This study proposed the thermo-hydro-mechanical phase-field model based on the micromechanically derived strain energy degradation with a new approach to update the porosity and a modified fixed stress split scheme with the thermal stress.
Furthermore, our results demonstrated the isotropic diffusion method's effectiveness in stabilizing the advection-dominated heat flux in hydraulic fracture.
Our model was verified against hydro-mechanical and thermo-hydro-mechanical problems and plane-strain hydraulic fracture propagation.
For hydraulic fracture propagation, we compared the results from our newly proposed porosity formulation against the existing model. 
We found that the new porosity formulation more accurately represents the sharp fracture behavior than the existing one.

Then, we illustrated the impacts of thermal stress on hydraulic fracture propagation through the examples with a single fracture and its interactions with a pre-existing weak interface under various temperature differences.
The higher the temperature difference, the less the critical pressure and the more fracture growth in both length and width. 
Moreover, because of this decrease in the critical pressure, higher temperature differences can suppress the hydro-mechanical attraction of the single fracture to the weak interface. 

As a future study, the thermal equilibrium assumption between the fracture and matrix may need to be relaxed as in~\cite{gelet2012thermo, gelet2013thermal}, and some analytical approaches may be implemented to ensure the continuity on the fracture boundary~\cite{noii2019phase, tran2013predicting}. 
Furthermore, the model may be extended to account for various working fluids such as supercritical CO$_2$~\cite{middleton2015shale} or liquid nitrogen~\cite{huang2020review}.
These unconventional fracturing fluids will likely go through the phase transition at the crack tip~\cite{yan2019staged}, which remains to be addressed in phase-field modelling of hydraulic fractures.

%% file: appendix.tex
\appendix

\section{Derivative of the Biot's modulus with respect to the phase field}
\label{A}
Recalling the incremental content of the fluid in Eq.~\ref{eq:simi_incremental_content_of_fluid}, Biot's modulus can also be written as
\begin{equation}
    \label{eq:Biot_modulus_new}
    \frac{1}{M_p} = \frac{(\zeta - \alpha \mrm{Tr}(\mbfs{\eps_e}))}{p}.
\end{equation}
Considering $p$ and $\zeta$ as the independent variables, the derivative of Biot's modulus with respect to the phase field is
\begin{equation}
    \label{eq:derivative_of_Biot_modulus_to_phi}
    \frac{\partial (1/M_p)}{\partial\upsilon} = \frac{-\mrm{Tr}(\mbfs{\eps_e})}{p}\frac{\partial\alpha}{\partial\upsilon},
\end{equation}
and the derivative of Biot's coefficient with respect to phase field writes
\begin{equation}
    \label{eq:derivative_of_Biot_coe_to_phi}
    \frac{\partial\alpha}{\partial\upsilon} = -2\upsilon(1-k)H(\mrm{Tr}\left(\varepsilon\right))\left(1 - \alpha_{\mathrm{m}}\right).
\end{equation}
Substituting Eq.~\ref{eq:derivative_of_Biot_coe_to_phi} into Eq.~\ref{eq:derivative_of_Biot_modulus_to_phi} yields
\begin{equation}
    \label{de_pf_Biot_modu_appendix}
    \frac{\partial 1/M_p}{\partial\upsilon} = \frac{2\eps_\mathrm{vol}}{p}\upsilon(1-k)H(\mrm{Tr}\left(\varepsilon\right))\left(1 - \alpha_{\mathrm{m}}\right)
    .
\end{equation}

\section{Discretization of THM phase-field modelling with FEM}
\label{B}
The variables $\upsilon$, $T$, $p$ and $\mbf{u}$ are defined at integration points as nodal values so that the field discretization form for the variables themselves and the respective gradients can be written as
\begin{equation}
		\begin{aligned}
			&\mbf{u}=\sum_{i=1}^{n}N_{i}^{u}u_i=\bm{N}_{u}\mbf{\hat{u}},\quad \upsilon=\sum_{i=1}^{n}N_{i}^{\upsilon}\upsilon_i=\bm{N}_{\upsilon}\hat{\upsilon},\quad p=\sum_{i=1}^{n}N_{i}^{p}p_i=\bm{N}_{p}\hat{p},\quad T=\sum_{i=1}^{n}N_iT_i=\bm{N}_{T}\hat{T}\\
			&\mbf{\eps}=\sum_{i=1}^{n}B_i^{u}u_i=\bm{B}_{u}{\hat{\mbf{u}}},\: \nabla \upsilon=\sum_{i=1}^{n}B_{i}^{\upsilon} \upsilon_{i}=\bm{B}_{\upsilon}\hat{\upsilon},\:\nabla p=\sum_{i=1}^{n}B_{i}^{p} p_{i}=\bm{B}_{p}\hat{p},\: \nabla T=\sum_{i=1}^{n}B_{i}^{T} T_{i}=\bm{B}_{T}\hat{T}\\
 	\end{aligned}
		\label{shap_f_varia}
	\end{equation}
where $\hat{T}$, $\hat{p}$, $\hat{\mbf{u}}$ and $\hat{\upsilon}$ represent the vectors of the integration point values in one element. 
Note that $\mbf{\eps}$ is a vector of independent strain variables but not a tensor, e.g., $\mbf{\eps}=\left[\eps_{xx},\:\eps_{yy},\:\eps_{xy}\right]^T$ in 2D problem. 
Shape functions $\bm{N}_u$, $\bm{N}_T$, $\bm{N}_p$ and $\bm{N}_\upsilon$ are represented as matrixes for vector field $\mbf{u}$ and a vector for scalar field like $p$, $\upsilon$ and $T$ and they are also used as test functions for each process:
\begin{align}
\begin{dcases}
        \bm{N}_u=\left[\begin{array}{cccccccc} N^u_1&0&...&N^u_i&0&...&N^u_n&0\\0&N^u_1&...&0&N^u_i&...&0&N^u_n
			\end{array}\right] \\
		\bm{N}_T=\left[\begin{array}{ccccc} N^T_1&...&N^T_i&...&N^T_n
			\end{array}\right] \\
        \bm{N}_p=\left[\begin{array}{ccccc} N^p_1&...&N^p_i&...&N^p_n
			\end{array}\right] \\
        \bm{N}_\upsilon=\left[\begin{array}{ccccc} N^\upsilon_1&...&N^\upsilon_i&...&N^\upsilon_n
			\end{array}\right] \\
\end{dcases}
.
\end{align}
The Galerkin finite element method considers the same shape functions in the present work, i.e., $N^u_i = N^T_i = N^p_i = N^\upsilon_i$. 
\begin{equation}
    \begin{aligned}
        \{\bm{T}\}_{m}=\{\bm{T}\}_{m-1}-\left[\bm{K^{TT}}\right]_{m-1}^{-1}\left\{\bm{r}^T\right\}_{m-1},
    \end{aligned}
    \label{N_P_T}
\end{equation}
\begin{equation}
    \begin{aligned}
        \{\bm{p}\}_{m}=\{\bm{p}\}_{m-1}-\left[\bm{K^{pp}}\right]_{m-1}^{-1}\left\{\bm{r}^p\right\}_{m-1},
    \end{aligned}
    \label{N_P_p}
\end{equation}
\begin{equation}
    \begin{aligned}
        \left\{\begin{array}{c}\mbf{u}\\ \bm{\upsilon} \end{array}\right\}_{m}=\left\{\begin{array}{c}\mbf{u}\\ \bm{\upsilon} \end{array}\right\}_{m-1}-\left[\begin{array}{cc}\bm{K}^{\mbf{u}\mbf{u}}&0\\ 0&\bm{K}^{\upsilon\upsilon} \end{array}\right]_{m-1}^{-1}\left\{\begin{array}{c}\bm{r}^{\mbf{u}}\\ \bm{r}^{\upsilon} \end{array}\right\}_{m-1}.
    \end{aligned}
    \label{NT}
\end{equation}

The residuals of each field are given as
\begin{flalign}
    &\bm{r}^T_{m-1}=\int_{\mit\Omega}\bm{N}_{T}^T\bm{N}_T(\rho c)_m\frac{\hat{T}^{k,m}-\hat{T}^{k-1}}{\Delta t} dV
    -\int_{\mit\Omega}\bm{B}_{T}^T\bm{N}_{T}\rho_{f}\hat{\mbf{q}_f}^{k,m-1} c_f T^{k,m}dV \notag \\
    &+\int_{\mit\Omega}\bm{B}_{T}^T\bm{B}_{T}\bm{\lambda}_\mathrm{eff} \hat{T}^{k,m}dV -\int_{\mit\Omega}\bm{N}^T_T Q_TdV+\int_{\partial_{N}\mit\Omega}\bm{N}^T_T \mbf{q}_{Tn}dS&
\end{flalign}
\begin{flalign}
        &\bm{r}^p_{m-1}=\int_{\mit\Omega} \alpha\bm{N}_{p}^T \frac{\varepsilon_v(\hat{\bm{u}}^{k,m-1})-\varepsilon_v(\hat{\bm{u}}^{k-1})}{\Delta t}dV + \int_{\mit\Omega} \bm{N}_{p}^T\frac{1}{M_p}\bm{N}_{p} dV \frac{\hat{p}^{k,m}-\hat{p}^{k-1}}{\Delta t} \notag \\
    &+\int_{\mit\Omega} \bm{N}_{p}^T\frac{1}{M_T}\bm{N}_{T}\frac{\hat{T}^{k,m}-\hat{T}^{k-1}}{\Delta t}dV +\int_{\mit\Omega}\bm{B}_{p}^T\frac{\mbf{K}}{\mu}\bm{B}_{p} dV \hat{p}^{k,m}+\int_{\mit\Omega} \frac{\alpha^2}{K}\bm{N}^T_p\bm{N}_pdV\frac{\hat{p}^{k,m}-\hat{p}^{k,m-1}}{\Delta t} \notag \\
    &+\int_{\mit\Omega} 3\alpha\alpha_s\bm{N}^T_p\bm{N}_T\frac{\hat{T}^{k,m}-\hat{T}^{k,m-1}}{\Delta t}dV-\int_{\mit\Omega}\bm{N}^T_p Q_fdV +\int_{\partial_{N}\mit\Omega}\bm{N}^T_pq_n dS&
\end{flalign}
\begin{flalign}
        \bm{r}^{\bm{u}}=&
        \int_{\mit\Omega}g(\upsilon)\bm{B}_{u}^T\mbf{\sigma}_e^{k,m}dV-\int_{\mit\Omega}\bm{B}_{u}^T\alpha(\upsilon)\bm{N}_{p}\hat{p}^{k,m}\bm{I}dV -\int_{\mathcal{C}_{N}}\bm{N}^T_{\mbf u}\bar{\mbf{t}} dS&
\end{flalign}
\begin{flalign}
        \bm{r}^{\upsilon}=&
        \int_{\mit\Omega}2(1-k) \bm{N}_{\upsilon}^T\bm{N}_{\upsilon}\psi_+(\mbf{u})^{k,m-1}\,dV\hat{\upsilon}+\frac{G_{c}}{4c_{n}}\int_{\mit\Omega}[-\frac{n}{\ell}(1-\bm{N}_{\upsilon}\hat{\upsilon})^{n-1}\bm{N}_{\upsilon}^T+2\ell\bm{B}_{\upsilon}^T\bm{B}_{\upsilon}\hat{\upsilon}]dV \notag \\
        +&\int_{\mit\Omega}2(1-k)\bm{N}_{\upsilon}^T\frac{(\bm{N}_{p}\hat{p})^2}{2}H(Tr\left(\varepsilon\right))\left(\frac{\alpha_m}{K_s}\right)dV\hat{\upsilon}&
\end{flalign}
and the Jacobian matrices for each field are given as
\begin{flalign}
    &\bm{K}^{TT}=\frac{\partial \bm{r}^{T}}{\partial \hat{T}}=\int_{\mit\Omega}\bm{N}_{T}^T\bm{N}_T(\rho c)_m\frac{1}{\Delta t} dV + \int_{\mit\Omega}\bm{B}_{T}^T\bm{\lambda}_\mathrm{eff}\bm{B}_{T}dV + \int_{\mit\Omega}\bm{B}_{T}^T\bm{N}_{T}\rho_{f}\hat{\mbf{q}_f}^{k,m-1} c_f \, dV&
\end{flalign}
\begin{flalign}
    &\bm{K}^{pp}=\frac{\partial \bm{r}^{p}}{\partial \hat{p}}=\int_{\mit\Omega} \bm{N}_{p}^T\frac{1}{M_p}\bm{N}_{p} dV \frac{1}{\Delta t}+\int_{\mit\Omega} \frac{\alpha^2}{K}\bm{N}^T\bm{N}dV\frac{1}{\Delta t}+\int_{\mit\Omega}\bm{B}_{p}^T\frac{\bm{K}}{\mu}\bm{B}_{p} dV&
\end{flalign}
\begin{flalign}
    &\bm{K}^{\bm{u}\bm{u}}=\frac{\partial \bm{r}^{\bm{u}}}{\partial \hat{\bm{u}}}=\int_{\mit\Omega}g(\upsilon)\bm{B}_{u}^T\mathbb{C}_+\bm{B}_{u}dV&
\end{flalign}
\begin{flalign}
    &\bm{K}^{\upsilon\upsilon}=\frac{\partial \bm{r}^{\upsilon}}{\partial \hat{\upsilon}}=\int_{\mit\Omega}2(1-k)\psi_+(\mbf{u})^{k,m-1}\bm{N}_{\upsilon}^T\bm{N}_{\upsilon}dV
        +\int_{\mit\Omega}2(1-k)\bm{N}_{\upsilon}^T\frac{(\bm{N}_{p}\hat{p})^2}{2}H(\Tr{\eps})\left(\frac{\alpha_m}{K_s}\right)dV\\
        &+\frac{G_{c}}{4c_{n}}\int_{\mit\Omega}[\frac{n(n-1)}{\ell}(1-\bm{N}_{\upsilon}\hat{\upsilon})^{n-2}\bm{N}_{\upsilon}^T\bm{N}_{\upsilon}+2\ell \bm{B}_{\upsilon}^T \bm{B}_{\upsilon}]dV&
\end{flalign}

\section{Approximation of the fracture width}
\label{C}
The fracture width expression of Eq.~\eqref{eq:frac_width0} is accurate for hydraulic fracturing under hydro-mechanical coupling because, without the thermal effect, $\eps_1$ can be approximated closely by the volumetric strain, i.e., $\eps_x + \eps_y$ in a 2D problem~\cite{YOU2023116305}. 
However, the thermal strain is isotropic and deforms in both the crack normal and tangential directions.
Thus, $\eps_1$ deviates from $\eps_\mrm{vol}$ when the thermal strain is involved.

\begin{figure}[H]
		\centerline{\includegraphics[scale=0.55]{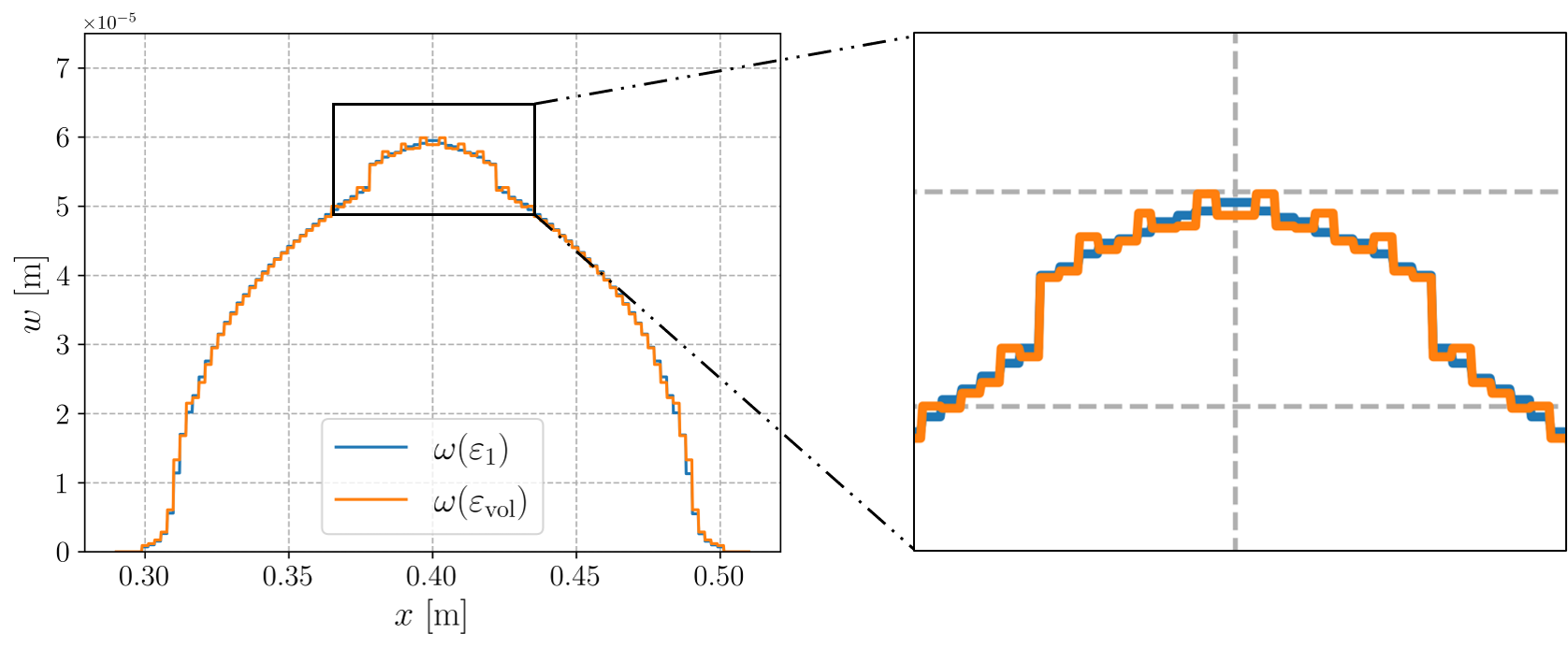}}
    \caption{\centering{Approximation of the fracture width with $\eps_1$ and $\eps_\mrm{vol}$.}}
		\label{fig:s_30K_width}
\end{figure}

To demonstrate this point, we simulated the single fracture model presented in section~\ref{sec:single-cracked model} with a temperature difference of 30~K using both Eq.~\eqref{eq:frac_width0} and Eq.~\eqref{eq:frac_width} to compute the fracture width (Fig.~\ref{fig:s_30K_width}).
The results show that the fracture width calculated by Eq.~\eqref{eq:frac_width0} exhibits an oscillating fracture width profile around the injection point.
